\pdfoutput=1
\newif\ifpersonal
\documentclass[11pt,a4paper,dvipsnames,x11names]{amsart}

\usepackage{amsmath,amsthm,amssymb,mathrsfs,mathtools,tensor,thmtools} 
\usepackage[mathcal]{eucal}
\usepackage[all,cmtip]{xy}

\usepackage[LGR,T1]{fontenc} 
\usepackage[utf8]{inputenc} 
\usepackage{CJKutf8}

\usepackage{microtype,inconsolata} 
\usepackage{enumerate,comment,braket,xspace,tikz,tikz-cd,csquotes} 
\usetikzlibrary{shapes.geometric}
\usepackage[centering,vscale=0.7,hscale=0.7]{geometry}
\usepackage{appendix}
\usepackage[hidelinks,colorlinks=true,linkcolor=Red4,citecolor=RoyalBlue,hypertexnames=false]{hyperref}

\usepackage[capitalize]{cleveref}
\usepackage{xcolor}
\usepackage{xpatch}

\usepackage{mathpazo}
\usepackage{euler}

\usepackage{stackrel}
\usepackage{faktor}
\usetikzlibrary{decorations.markings}

\numberwithin{equation}{subsection}
\theoremstyle{plain}
\newtheorem{theorem}[equation]{Theorem}
\newtheorem{lemma}[equation]{Lemma}

\newtheorem*{example*}{Example}
\newtheorem*{definition*}{Definition}
\newtheorem{proposition}[equation]{Proposition}

\newtheorem{corollary}[equation]{Corollary}
\newtheorem{assumption}[equation]{Assumption}

\theoremstyle{definition}
\newtheorem{definition}[equation]{Definition}
\newtheorem{notation}[equation]{Notation}
\newtheorem{example}[equation]{Example}
\newtheorem{remark}[equation]{Remark}

\newtheorem{recollection}[equation]{Recollection}

\newtheorem{construction}[equation]{Construction}

\ifpersonal
\newcommand{\personal}[1]{\textcolor[rgb]{0,0,1}{(Personal: #1)}}
\newcommand{\discussion}[1]{\textcolor{violet}{(Discussion: #1)}}
\else
\newcommand{\personal}[1]{\ignorespaces}
\newcommand{\discussion}[1]{\ignorespaces}
\fi

\usetikzlibrary{decorations.markings}
\tikzset{double line with arrow/.style args={#1,#2}{decorate,decoration={markings,%
mark=at position 0 with {\coordinate (ta-base-1) at (0,1pt);
\coordinate (ta-base-2) at (0,-1pt);},
mark=at position 1 with {\draw[#1] (ta-base-1) -- (0,1pt);
\draw[#2] (ta-base-2) -- (0,-1pt);
}}}}

\def\G{\mathbb{G}}
\def\N{\mathbb{N}}

\def\ZZ{\mathbb{Z}}

\def\C{\mathcal{C}}
\def\D{\mathcal{D}}
\def\E{\mathcal{E}}

\def\X{\mathcal{X}}
\def\Y{\mathcal{Y}}
\def\Z{\mathcal{Z}}

\def\M{\mathcal{M}}
\def\OO{\mathcal{O}}

\def\Spl{\mathbf{Spl}}

\def\Prlo{\mathrm{Pr}^{\mathrm{L}, \omega}}

\def\Funst{\mathrm{Fun}^{\mathrm{st}}}
\def\FunR{\mathrm{Fun}^{\mathrm{R}}}

\DeclareMathOperator{\QCoh}{QCoh}
\DeclareMathOperator{\dAff}{dAff}
\DeclareMathOperator{\dSt}{dSt}

\DeclareMathOperator{\Aff}{Aff}
\DeclareMathOperator{\fib}{fib}
\DeclareMathOperator{\cof}{cof}

\DeclareMathOperator{\forget}{forget}

\DeclareMathOperator{\Hom}{Hom}

\DeclareMathOperator{\Fun}{Fun}

\DeclareMathOperator{\HH}{HH}

\DeclareMathOperator{\Mod}{Mod}

\DeclareMathOperator{\Ext}{Ext}

\DeclareMathOperator{\St}{St}

\DeclareMathOperator{\Map}{Map}

\DeclareMathOperator{\Cat}{\mathrm{Cat}}
\DeclareMathOperator{\Spc}{Spc}
\DeclareMathOperator{\ev}{ev}
\DeclareMathOperator{\coev}{coev}
\DeclareMathOperator{\Id}{Id}

\DeclareMathOperator{\Perf}{Perf}
\DeclareMathOperator{\APerf}{APerf}

\DeclareMathOperator{\DCoh}{DCoh}
\DeclareMathOperator{\IndCoh}{IndCoh}

\DeclareMathOperator{\Spec}{Spec}

\DeclareMathOperator{\simple}{sim}
\DeclareMathOperator{\spl}{spl}

\newcommand{\colim}{\operatornamewithlimits{colim}}

\makeatletter
\renewenvironment{proof}[1][\relax]{\par
  \pushQED{\qed}%
  \normalfont \topsep6\p@\@plus6\p@\relax
  \trivlist
  \item[\hskip\labelsep\itshape
    \ifx#1\relax \proofname\else\proofname{} of #1\fi\@addpunct{.}]\ignorespaces
}{%
  \popQED\endtrivlist\@endpefalse
}
\makeatother

\makeatletter
\let\save@mathaccent\mathaccent
\newcommand*\if@single[3]{%
	\setbox0\hbox{${\mathaccent"0362{#1}}^H$}%
	\setbox2\hbox{${\mathaccent"0362{\kern0pt#1}}^H$}%
	\ifdim\ht0=\ht2 #3\else #2\fi
}
\newcommand*\rel@kern[1]{\kern#1\dimexpr\macc@kerna}
\newcommand*\widebar[1]{\@ifnextchar^{{\wide@bar{#1}{0}}}{\wide@bar{#1}{1}}}
\newcommand*\wide@bar[2]{\if@single{#1}{\wide@bar@{#1}{#2}{1}}{\wide@bar@{#1}{#2}{2}}}
\newcommand*\wide@bar@[3]{%
	\begingroup
	\def\mathaccent##1##2{%
		\let\mathaccent\save@mathaccent
		\if#32 \let\macc@nucleus\first@char \fi
		\setbox\z@\hbox{$\macc@style{\macc@nucleus}_{}$}%
		\setbox\tw@\hbox{$\macc@style{\macc@nucleus}{}_{}$}%
		\dimen@\wd\tw@
		\advance\dimen@-\wd\z@
		\divide\dimen@ 3
		\@tempdima\wd\tw@
		\advance\@tempdima-\scriptspace
		\divide\@tempdima 10
		\advance\dimen@-\@tempdima
		\ifdim\dimen@>\z@ \dimen@0pt\fi
		\rel@kern{0.6}\kern-\dimen@
		\if#31
		\overline{\rel@kern{-0.6}\kern\dimen@\macc@nucleus\rel@kern{0.4}\kern\dimen@}%
		\advance\dimen@0.4\dimexpr\macc@kerna
		\let\final@kern#2%
		\ifdim\dimen@<\z@ \let\final@kern1\fi
		\if\final@kern1 \kern-\dimen@\fi
		\else
		\overline{\rel@kern{-0.6}\kern\dimen@#1}%
		\fi
	}%
	\macc@depth\@ne
	\let\math@bgroup\@empty \let\math@egroup\macc@set@skewchar
	\mathsurround\z@ \frozen@everymath{\mathgroup\macc@group\relax}%
	\macc@set@skewchar\relax
	\let\mathaccentV\macc@nested@a
	\if#31
	\macc@nested@a\relax111{#1}%
	\else
	\def\gobble@till@marker##1\endmarker{}%
	\futurelet\first@char\gobble@till@marker#1\endmarker
	\ifcat\noexpand\first@char A\else
	\def\first@char{}%
	\fi
	\macc@nested@a\relax111{\first@char}%
	\fi
	\endgroup
}
\makeatother

\usetikzlibrary{decorations.markings} 
\tikzset{
  closed/.style = {decoration = {markings, mark = at position 0.5 with { \node[transform shape, xscale = .8, yscale=.4] {/}; } }, postaction = {decorate} },
  open/.style = {decoration = {markings, mark = at position 0.5 with { \node[transform shape, scale = .7] {$\circ$}; } }, postaction = {decorate} }
}

\makeatletter
\tikzcdset{
  open/.code     = {\tikzcdset{hook, circled};},
  closed/.code   = {\tikzcdset{hook, slashed};},
  open'/.code    = {\tikzcdset{hook', circled};},
  closed'/.code  = {\tikzcdset{hook', slashed};},
  circled/.code  = {\tikzcdset{markwith = {\draw (0,0) circle (.375ex);}};},
  slashed/.code  = {\tikzcdset{markwith = {\draw[-] (-.4ex,-.4ex) -- (.4ex,.4ex);}};},
  markwith/.code ={
    \pgfutil@ifundefined%
    {tikz@library@decorations.markings@loaded}%
    {\pgfutil@packageerror{tikz-cd}{You need to say %
      \string\usetikzlibrary{decorations.markings} to use arrows with markings}{}}{}%
    \pgfkeysalso{/tikz/postaction = {
      /tikz/decorate,
      /tikz/decoration={markings, mark = at position 0.5 with {#1}}}
    }
  },
}
\makeatother

\newcommand{\leftrarrows}{\mathrel{\raise.75ex\hbox{\oalign{%
  $\scriptstyle\leftarrow$\cr
  \vrule width0pt height.5ex$\hfil\scriptstyle\relbar$\cr}}}}
\newcommand{\lrightarrows}{\mathrel{\raise.75ex\hbox{\oalign{%
  $\scriptstyle\relbar$\hfil\cr
  $\scriptstyle\vrule width0pt height.5ex\smash\rightarrow$\cr}}}}
\newcommand{\Rrelbar}{\mathrel{\raise.75ex\hbox{\oalign{%
  $\scriptstyle\relbar$\cr
  \vrule width0pt height.5ex$\scriptstyle\relbar$}}}}

\makeatletter
\def\leftrightarrowsfill@{\arrowfill@\leftrarrows\Rrelbar\lrightarrows}
\newcommand{\xleftrightarrows}[2][]{\ext@arrow 3399\leftrightarrowsfill@{#1}{#2}}
\makeatother

\NewCommandCopy{\notocsection}{\section}
\xpatchcmd{\notocsection}{{1}}{{1001}}{}{}

\begin{document}

\title{Shifted Lagrangian structures and syzygies on Calabi--Yau varieties}
\author{Enrico Lampetti and Luca Morstabilini}
\address{ENRICO LAMPETTI, UNIVERSITY OF SOUTHERN CALIFORNIA, LOS ANGELES, CA, USA}
\email{lampetti@usc.edu}

\address{LUCA MORSTABILINI, SISSA, TRIESTE, ITALY}
\email{lmorstab@sissa.it}

\subjclass[2020]{}

\begin{abstract}
	In \cite{FMR25}, Fantechi--Miró-Roig show how to construct new Lagrangian subspaces of the moduli of simple perfect sheaves on a K3 surface from a given one.
	We reinterpret their construction using the tools and flexibility of derived algebraic geometry.
	This allows us to generalize from K3 surfaces to (possibly singular) Calabi--Yau schemes of arbitrary dimensions.
\end{abstract}

\maketitle

\tableofcontents

\section{Introduction}
\addtocontents{toc}{\protect\setcounter{tocdepth}{1}}

    In his seminal paper \cite{Muk84}, Mukai proved that simple coherent sheaves on a K3 or abelian variety $X$ organize into a smooth symplectic algebraic space $\underline{Spl}(X)$.
    Later, Inaba \cite{Ina11} generalized Mukai's Theorem by showing that the algebraic space $\underline{Perf}^{sim}(X)$ parametrizing simple complexes with bounded coherent cohomology on $X$ is smooth and symplectic.\medskip

    The advent of homotopical technique in algebraic geometry, developed by Toën--Vezzosi \cite{TV05, TV08} and Lurie \cite{Lur17, Lur18}, allowed to vastly generalize the results of Mukai \cite{Muk84} and Inaba \cite{Ina11}.
    Indeed, in the breakthrough paper \cite{PTVV13}, Pantev--Toën--Vezzosi--Vaquié introduced the notion of shifted symplectic structure, and showed that the \textit{derived stack} $\mathbf{Perf}(X)$ of perfect complexes on a smooth and proper scheme $X$ with trivial canonical bundle admits a $(2- \dim X)$-shifted symplectic form.
    Examples of shifted Lagrangian substacks of $\mathbf{Perf}(X)$ were given by Brav--Dyckerhoff \cite{BD19, BD21}.
    Notice that for $X$ smooth, complexes with bounded coherent cohomology and perfect complexes agree, but this is not the case in general and will be important later on.
    Furthermore, if $\dim X = 2$, it was also shown in \cite{PTVV13} that one could recover the existence of a classical symplectic structure on the classical algebraic space $\underline{Perf}^{sim}(X)$ from the existence of a $0$-shifted symplectic structure on $\mathbf{Perf}(X)$, hence reproving the aforementioned results of Mukai \cite{Muk84} and Inaba \cite{Ina11}.
    One of the main advantages of derived algebraic geometry is that it gives a framework where checking the closedness of the symplectic forms above does not encompass heavy computations (see the discussion before \cref{recollection:gerbes_simple} and the subsequent discussion).
\medskip

    Going back to classical algebraic geometry, Fantechi--Miró-Roig \cite{FMR25, FMR25b} studied Lagrangian sub-algebraic spaces of $\underline{Perf}^{sim}(X)$, with $X$ a K3 surface.
    More precisely, in \cite{FMR25} a way of constructing new Lagrangian sub-algebraic spaces from a fixed one is given, and is related to moduli spaces of syzygy bundles \cite{FMR24}.
    Moreover, the paper \cite{FMR25b} extends the result of \cite{FMR25} to possibly singular K3 surfaces.
    In \cite{FMR25b}, the authors need to impose a smoothability condition in order for the argument to go through.
    The main goal of this paper is to extend their results in two different directions:
    \begin{enumerate}
        \item we consider Lagrangians not in the algebraic space $\underline{Spl}^(X)$ of simple perfect sheaves, but in the whole stack $\mathbf{Perf}(X)$ of perfect complexes;
        \item we allow $X$ to be \emph{any} singular K3 surface.
        More generally, we allow $X$ to be any \textit{possibly singular} Calabi--Yau scheme of \textit{arbitrary dimension}.
    \end{enumerate}
    Our results hold for \emph{strong Calabi--Yau schemes}, that is, schemes $X$ with trivial canonical bundle and such that $\mathrm{H}^i(X, \OO_X) \simeq 0$ for $i \neq 0, \dim X$.
    We refer to \cref{def:strong_CY_scheme} for a precise definition of these objects.
    In dimension 2, they are precisely (singular) K3 surfaces and we are moreover able to recover the results of \cite{FMR25b}.
    \medskip
    
    Our approach makes extensive use of derived algebraic geometry and higher category theory. The starting point is the construction of a derived enhancement of the stack of syzygy bundles considered in \cite{FMR24}.
    For this, let us fix a proper scheme $X$ over an algebraically closed field $k$ of characteristic zero.
    Consider the derived stack $\mathbf{dSyz}(X)$, whose $k$-points correspond to simple perfect sheaves $F \in \QCoh(X_A)^\heartsuit$ and vector subspace $W \subset \mathrm{H}^0(X,F)$ that globally generate $F$, see \cref{rem:comparison_dSyz}.
    Actually, we will consider a bigger derived stack $\widetilde{\mathbf{Syz}}(X)$, of which $\mathbf{dSyz}(X)$ is an open substack.
    The $A$-points of $\widetilde{\mathbf{Syz}}(X)$ correspond to morphisms $g \colon p_A^\ast W \to F$, where $F \in \Perf(X_A)$, $W \in \Perf(A)$.
    Here, we let $X_A \coloneqq X \times_k \Spec(A)$ and $p_A \colon X_A \to \Spec(A)$ the projection.
    The derived stack $\widetilde{\mathbf{Syz}}(X)$ is equipped with two maps $\ev_1, \fib \colon \widetilde{\mathbf{Syz}}(X) \to \mathbf{Perf}(X)$, defined respectively by
\[
(g \colon p_A^\ast W \to F) \xmapsto{\ev_1} F, \qquad (g \colon p_A^\ast W \to F) \xmapsto{\fib} \fib(g) \ .
\]
    The following is our first main result.

\begin{theorem}[{\cref{thm:main_DCoh}}]\label{thm:main_DCoh_intro}
    Let $X$ be a strong Calabi--Yau scheme of dimension $d > 0$ (\cref{def:strong_CY_scheme}) and $i \colon \mathbf{\mathcal{L}} \to \mathbf{Perf}(X)$ be a $(1-d)$-shifted Lagrangian morphism.
    Consider the diagram in $\dSt_k$ with pullback square
\[
\begin{tikzcd}[sep=small]
	{\widetilde{\mathbf{Syz}}^{\mathcal{L}}(X)} && {\widetilde{\mathbf{Syz}}(X)} && {\mathbf{Perf}(X)}
	\\
    \\
	{\mathbf{\mathcal{L}}} && {\mathbf{Perf}(X)} 
    \arrow["\fib", from=1-3, to=1-5]
	\arrow["q", from=1-1, to=1-3]
	\arrow[from=1-1, to=3-1]
	\arrow["\ev_1"', from=1-3, to=3-3]
	\arrow["i"', from=3-1, to=3-3]
\end{tikzcd}
\]
	Then, the morphism of derived stacks
\[
\fib \circ q \colon \widetilde{\mathbf{Syz}}^{\mathbf{\mathcal{L}}}(X) \to \mathbf{Perf}(X)
\]
is $(1-d)$-shifted Lagrangian.
\end{theorem}

    \cref{thm:main_DCoh_intro} follows from the composition of Lagrangian spans \cite[Theorem 4.4]{Cal15} obtained by Calaque, together with the following result:

\begin{theorem}[{\cref{coro:Lagr_strong}}]\label{coro:Lagr_strong_intro}
    Let $X$ be a proper strong Calabi--Yau scheme of dimension $d > 0$.
	Then, the morphism in $\dSt_k$
\[
p^\ast \colon \mathbf{Perf}(k) \to \mathbf{Perf}(X) 
\]
carries a $(1-d)$-shifted Lagrangian structure.
\end{theorem}

    Notice that \cref{coro:Lagr_strong_intro} (as well as the composition of Lagrangian spans \cite[Theorem 4.4]{Cal15}) is a fully derived statement, even for $\dim X = 2$.\medskip
    
    One can prove \cref{coro:Lagr_strong} directly, as shown in Appendix \ref{appendix:the_case_of_Perf}.
    Nevertheless, we prove a finer statement in \cref{prop:CY_strong}.
    More precisely, we prove the existence of a non-commutative shifted Lagrangian  morphism, that is, a relative left Calabi--Yau morphism \cite{BD19} on categories of Ind-coherent sheaves.
    Then, appealing to the main result of Brav--Dyckerhoff \cite{BD21}, we get a shifted Lagrangian morphism by passing to the derived moduli of objects, as introduced by Toën-Vaquié \cite{TV07}.
    This strategy relies on \cref{thm:comparison_moduli_Perf_DCoh}, which might be of independent interest, exhibiting a canonical equivalence between the moduli of objects of Ind-Coherent sheaves and the derived stack of perfect complexes.
    A variant for non-proper schemes is also proven. 
    More precisely, we prove the following. 
    \begin{theorem}[{\cref{thm:comparison_moduli_Perf_DCoh}}]\label{thm:moduli_of_objects_intro}
    Let $X$ be a classical separated scheme, locally of finite presentation over $k$.
    Then, the moduli of objects $\M_{\IndCoh(X)}$ of Ind-coherent sheaves on $X$ is equivalent to the derived stack $\mathbf{Perf}_{prop}(X)$ parametrizing perfect complexes on $X$ with proper support.
    \end{theorem}
    \cref{thm:moduli_of_objects_intro} seems to be a result expected by experts, and is implicitly used in the introduction of \cite{BD21}. We could not locate a proof in the literature, so we provide a full proof in Appendix \ref{Appendix:equivalence}.
    Our proof relies on the theory of integral transforms for coherent sheaves of Ben-Zvi--Nadler--Preygel \cite{BZNP17}.
    During the writing of this paper, a proof of \cref{thm:moduli_of_objects_intro} in the case where $X$ is smooth over a smooth base appeared in \cite{Huyhn26}.
    Notice that, in the setting of \cite{Huyhn26}, one has $\IndCoh(X) \simeq \QCoh(X)$, which substantially simplifies the proof.\medskip

    Finally, we show in \cref{section:comparison_with_FMR} how to deduce results for classical algebraic spaces of perfect sheaves from the results on derived stacks of \cref{thm:main_DCoh_intro}, hence generalizing \cite{FMR25b} by removing the aforementioned smoothability assumption.
    For this, we show that the \textit{derived symplectic rigidification}, introduced in \cite{Park26}, is actually a smooth underived stack and recovers the algebraic space studied in \cite{FMR25b}.
    This phenomenon is special to dimension 2: in higher dimensions, derived structure naturally appears in the quotient algebraic space.
    More precisely, in the setting of \cref{thm:main_DCoh_intro}, let $d = 2$ and $k = \mathbb{C}$.
    Consider the morphism
\[
\fib \circ q \colon \mathbf{dSyz}^{\mathcal{L}}(X) \times_{\mathbf{Perf}(X)} \Spl(X) \to \Spl(X)
\]
obtained by pulling back $\fib \circ q \colon \widetilde{\mathbf{Syz}}^{\mathbf{\mathcal{L}}}(X) \to \mathbf{Perf}(X)$ along the open immersion $\Spl(X)\subset \mathbf{Perf}(X)$ (\cref{def:simple_sheaves}) and by restricting the source along the open immersion $\mathbf{dSyz}(X) \subset \widetilde{\mathbf{Syz}}(X)$ (\cref{lem:openness_dSyz}).
	Then, we have:
\begin{theorem}[\cref{coro:classic_Lagr_from_derived_Lagr}]\label{intro:classic_Lagr_from_derived_Lagr}
	Assume that $\mathbf{\mathcal{L}}$ carries a $B\G_m$-action and that $i$ is $B\G_m$-equivariant.
	Then, the schematic image of the morphism
\[
t_0(\fib \circ q)^{rig} \colon t_0\big(\mathbf{dSyz}^{\mathcal{L}}(X) \times_{\mathbf{Perf}(X)} \Spl(X)\big) / B\G_m \to \underline{Spl}(X)
\]
is a classical Lagrangian sub-algebraic space, where $t_0\big(\mathbf{dSyz}^{\mathcal{L}}(X) \times_{\mathbf{Perf}(X)} \Spl(X)\big) / B\G_m$ is the classical rigidification of $t_0\big(\mathbf{dSyz}^{\mathcal{L}}(X) \times_{\mathbf{Perf}(X)} \Spl(X)\big)$.
\end{theorem}
	
	\cref{intro:classic_Lagr_from_derived_Lagr} is a first generalization of \cite[Theorem 1.3]{FMR25b}, where the input is a shifted Lagrangian morphism $\mathcal{L} \to \mathbf{Perf}(X)$.
	Since $\Spl(X) \subset \mathbf{Perf}(X)$ is an open immersion, one can also take as input a shifted Lagrangian morphism $\mathcal{L} \to \Spl(X)$.
	In \cref{thm:classical_results} we obtain the same result by starting with a classical Lagrangian sub-algebraic space $L \to \underline{Spl}(X)$, hence proving a faithful generalization of \cite[Theorem 1.3]{FMR25b}.
	See \cref{rem:FMR_comparison}.\medskip

    Our techniques allow us to simplify and strengthen the arguments of \cite{FMR25, FMR25b}, and should therefore help in the study of moduli spaces of sheaves on more general singular $K3$-surfaces. Moreover, the ability to construct large families of Lagrangians allows, through the derived Lagrangian intersection theorem \cite[Theorem 2.9]{PTVV13} to exhibit examples of $(-1)$-shifted symplectic stacks, i.e., spaces where DT invariants can be defined \cite{darbouxtheorem}.

\subsection{Linear overview}

    In \cref{section_moduli_of_objects}, we recall the framework and results on moduli of objects developed by Toën--Vaquié \cite{TV07}.
    In \cref{subsection:moduli_of_objects_Perf}, we specialize the discussion to the case of perfect complexes, which is the focus of the present paper.\\ \indent
    In \cref{section:CY_structures}, we recall the definition of relative left Calabi--Yau structure, following \cite{BD19, Chr23}.
    We recall in \cref{thm:CY_lagr} the main result of \cite{BD21}, stating that relative left Calabi--Yau structures are non-commutative analogues of shifted Lagrangian morphisms.
    In \cref{subsection:CY_algebraic_geometry}, we will specialize our discussion to categories of Ind-coherent sheaves on Calabi--Yau schemes.
    We will use the results therein to deduce statements about moduli of perfect complexes.\\ \indent
    In \cref{section:main_results}, we use the results proven in \cref{section:CY_structures} to deduce our main result about stacks of perfect complexes on strong Calabi--Yau schemes, namely, \cref{thm:main_DCoh}.\\ \indent
    In \cref{section:comparison_with_FMR}, we show how to recover and generalize the results of \cite{FMR25b} from \cref{thm:main_DCoh}.\\ \indent
    In Appendix \ref{Appendix_FM}, we collect results about Fourier--Mukai transforms.
    This leads to a proof of \cref{thm:comparison_Perf_corpus}, on which our paper relies.
    \cref{thm:comparison_Perf_corpus} allows us to identify the moduli of perfect complexes with proper support on a scheme $X$ with the moduli of objects of Ind-coherent sheaves on $X$.\\ \indent
    In Appendix \ref{appendix:the_case_of_Perf}, we give a fully commutative proof of \cref{coro:Lagr_strong}, that is, a proof that does not rely on the theory of Calabi--Yau structures.

\subsection*{Notations}
	We introduce the following running notations.
\begin{itemize}\itemsep=0.2cm
    \item We fix $k$ an algebraically closed field of characteristic $0$;
    \item All (derived) schemes are quasi-compact quasi-separated;
    \item $\mathrm{Grpd}$ is the category of groupoids;
    \item $\mathrm{Spc}$ is the $\infty$-category of spaces;
    \item $\mathrm{Sp}$ is the $\infty$-category of spectra;
    \item $\Cat_\infty$ is the $\infty$-category of $\infty$-categories;
    \item Given $\C \in \Cat_\infty$, we let $\C^\omega\subset \C$ be the full subcategory spanned by compact objects;
    \item $(-)^\simeq \colon \Cat_\infty \to \Spc$ is the functor that assigns to an $\infty$-category its maximal $\infty$-sub-groupoid;
    \item $\Pr^{\mathrm{L}}_k$ is the $\infty$-category of presentable $k$-linear $\infty$-categories with left adjoint functors;
    \item $\Prlo_k$ is the $\infty$-category of presentable compactly generated $k$-linear $\infty$-categories with left adjoint functors preserving compact objects;
    \item $\Aff_k$ is the opposite of the category of $k$-algebras;
    \item $\dAff_k$ is the opposite of the $\infty$-category of simplicial $k$-algebras;
    \item the category of stacks over $k$ is the full subcategory
    \[
    \mathrm{St}_k \subset \Fun(\Aff_k^{op}, \mathrm{Grpd})
    \]
    spanned by stacks over $k$ for the étale topology;
    \item the $\infty$-category of derived stacks over $k$ is the full subcategory
    \[
    \dSt_k \subset \Fun(\dAff_k^{op}, \mathrm{Spc})
    \]
    spanned by derived stacks over $k$ for the hypercomplete étale topology;
    \item denote by $\tau_{\leq 1} \colon \mathrm{Spc} \to \mathrm{Grpd}$ the functor that assigns to a space its homotopy category, which is a groupoid.
    Denote by $i \colon \Aff_k \to \dAff_k$ the canonical inclusion.
    The composition 
    \[
    \Fun(\dAff_k^{op}, \mathrm{Spc}) \xrightarrow{- \circ i} \Fun(\Aff_k^{op}, \mathrm{Spc}) \xrightarrow{\tau_{\leq 1} \circ -} \Fun(\Aff_k^{op}, \mathrm{Grpd})
    \]
    restricts to a functor $t_0 \colon \dSt_k \to \St_k$.
    We will refer to $t_0$ as the truncation functor.
    \item
	For $X \in \dSt_k$ and $\Spec(A) \in \dAff_k$, we denote
\[
X_A \coloneqq X \times_k \Spec(A) \ .
\]
\end{itemize}

\addtocontents{toc}{\protect\setcounter{tocdepth}{0}}
\section*{Acknowledgments}

    This paper would not exist if the authors had not attended the conference "Stacks and Moduli" at IIT Madras in January 2026, where Fantechi and Miró-Roig asked whether the smoothability condition in \cite{FMR25b} could be dropped.
    We thank the organizers for the beautiful conference and the IIT Madras for the hospitality during our stay.\\ \indent
    E.L. wishes to thank his advisors Mauro Porta and Jean-Baptiste Teyssier for their support during his thesis and while writing this article. 
    E.L. also thanks Valerio Melani and Gabriele Vezzosi for fruitful discussions about the content of this paper.
    L.M. wishes to thank his advisors Francesco Sala and Andrea Ricolfi for their support and for useful discussions.\\ \indent
    The authors wish to thank Mauro Porta for fruitful discussions about the content of Appendix \ref{Appendix:equivalence}.
    The authors also thank Matteo Verni for a fruitful discussion about \cref{example:Strong_CY_from_anticanonical_div}.
    They also wish to thank Hyeonjun Park and Jemin You for answering their questions about the content of \cref{section:shifted_sympl_rig} (see \cref{rem:PY_shifted_Lagr_rig}).

\section*{AI statement}
    Artificial intelligence was not used for the mathematics or the writing of this paper.
    AI was used for typo-checking, but not for proof checking.

\addtocontents{toc}{\protect\setcounter{tocdepth}{2}}

\section{Moduli of objects}\label{section_moduli_of_objects}

	In this Section we recall the main definition of moduli of objects and the main results needed in the rest of the paper.

\subsection{Families of objects}

	In this paragraph we recall properties of the tensor product of presentable $\infty$-categories, which is used to define families of objects. 
	We also discuss some finiteness properties one can impose on compactly generated presentable categories.

\begin{recollection}[{\cite[Section 4.8]{Lur17}}]\label{tensor_product}
	For $\C, \D \in \mathrm{Pr}^\mathrm{L}_k$, we can consider their tensor product
\[
\C \otimes_k \D \simeq \FunR_k(\C^{op}, \D)
\]
where the right hand side is the $\infty$-category of $k$-linear functors commuting with limits.
	This tensor product endows $\Pr^{\mathrm{L}}_k$ with a symmetric monoidal structure which restricts to $\Prlo_k$. 
	When $\C$ is compactly generated, there is a canonical equivalence
\[
\C \otimes_k \D \simeq \Funst_k((\C^{\omega})^{op}, \D)
\]
where the right hand side is the $\infty$-category of exact $k$-linear functors.
\end{recollection}

\begin{remark}\label{functoriality_tens_prod}
	Let $\C, \D, \E \in \Pr_k^{\mathrm{L}}$.
	Then an adjunction
\[
f \colon \D \leftrightarrows \E \colon g
\]
induces an adjunction 
\[
\Id_\C \otimes_k f \colon \C \otimes_k \D\leftrightarrows \C \otimes_k \E \colon \Id_\C \otimes_k g.
\]
	The functor $\Id_\C \otimes_k g$ corresponds to 
\[
g \circ - \colon \FunR_k(\C^{op}, \E) \to \FunR_k(\C^{op}, \D)
\]
under the equivalences of \cref{tensor_product}.
	Its left adjoint $\Id_\C \otimes_k f$ has no easy description in general.
	Notice indeed that post-composition with $f$ does not preserve functors that commute with limits.
	Nevertheless, if $\C$ is compactly generated one can see that the adjunction $\Id_\C \otimes_k f  \dashv \Id_\C \otimes_k g$ corresponds to 
\[
f \circ - \colon \Funst_k((\C^{\omega})^{op}, \D) \leftrightarrows \Funst_k((\C^{\omega})^{op}, \E) \colon g \circ -
\]
under the equivalences of \cref{tensor_product}.
\end{remark}

\begin{notation}
In the setting of \cref{functoriality_tens_prod}, we denote the functors $\Id_\C \otimes_k f, \Id_\C \otimes_k g$ respectively by $f_\C, g_\C$.
\end{notation}

	In the stable setting, the following result is a reformulation of \cite[Lemma 2.2.1]{SS03}.

\begin{lemma}[{\cite[Lemma 2.4]{AG14}}]\label{generators_stable}
	Let $\E \in \Prlo_k$ stable and $X \subset \E^\omega$ a set of compact objects.
	Then, the following are equivalent:
\begin{enumerate}\itemsep=0.2cm
	\item $\forall y \in \E$,
\[
\mathrm{Map}_\E(x,y) \simeq * \in \mathrm{Sp} \ \  \forall x \in X \Rightarrow y \simeq 0;
\]
	\item the $\infty$-category $\E^\omega$ is the smallest full subcategory of $\E$ containing $X$ and stable under finite colimits, shifts and retracts.
\end{enumerate}
\end{lemma}

\begin{definition}\label{compact_gen_def}
	Let $\E \in \Prlo_k$ stable and $X \subset \E^\omega$ a set of compact objects.
	We say that $X$ is a set of compact generators if it satisfies the equivalent conditions of \cref{generators_stable}.
\end{definition}

\begin{example}
	Let $A \in \mathrm{Alg}_k$.
	Then $A$ is a compact generator of $\Mod_A$ in the sense of \cref{compact_gen_def}.
\end{example}

\begin{definition}
	For $\C \in \mathrm{Pr}^\mathrm{L}_k$ and $X\in \dSt_k$, define the $\infty$-category of $X$-families of objects in $\C$ as
\[
\C_X \coloneqq  \C \otimes_k \QCoh(X) \ .
\]
\end{definition}

\begin{notation}
	Let $\C \in \Pr^{\mathrm{L}}_k$ and $S = \Spec (A) \in \dAff_k$.
	The $\infty$-category of $S$-families of objects in $\C$ will be denoted equivalently by $\C_S$ or by $\C_A$.
\end{notation}

\begin{definition}\label{pseudo_perfect}
	For $\C \in \Prlo_k$ and $X \in \dSt_k$, define the $\infty$-category of pseudo-perfect families over $X$ in $\C$ as the full subcategory
\[
\Funst((C^\omega)^{op}, \Perf(X)) \subset \C_X
\]
spanned by functors valued in $\Perf(X)$.
\end{definition}

	In order to introduce finiteness conditions on presentable categories, let us give the following:

\begin{lemma}[{\cite[Theorem D.7.0.7]{Lur18}}]\label{lemma:Prlo_dualizable}
	Let $\C \in \Prlo_k$.
	Then $\C$ is dualizable in $\Pr_k^{\mathrm{L}}$ with dual given by $\C^\vee \simeq \mathrm{Ind}((\C^\omega)^{op})$.
	In particular, up to a contractible space of choices, there are unique maps in $\Pr^{\mathrm{L}}_k$
\[
\mathrm{ev}_\C \colon  \C^\vee \otimes_k \C \to \Mod_k
\]
\[
\mathrm{coev}_\C  \colon  \Mod_k \to \C \otimes_k \C^\vee
\]
satisfying the triangular identities.
\end{lemma}

	We will need the following finiteness condition:

\begin{definition}\label{def_finiteness_conditions}
	Let $\C \in \Prlo_k$.
	We say that $\C$ is 
\begin{enumerate}\itemsep=0.2cm
    \item smooth if $\coev_\C$ preserves compact objects;
    \item proper if $\ev_\C$ preserves compact objects;
    \item of finite type if it is a compact object of $\Prlo_{k}$.
\end{enumerate}
\end{definition}

	The next result compares the different finiteness conditions of \cref{def_finiteness_conditions}.

\begin{proposition}[{\cite[Proposition II.2.6]{DPS25}}]\label{relation_finiteness_conditions}
	Let $\C \in \Prlo_k$.
	Then:
\begin{enumerate}\itemsep=0.2cm
    \item if $\C$ is of finite type, it is smooth;
    \item if $\C$ is smooth and proper, it is of finite type.
\end{enumerate}
Moreover, if $\C$ is smooth then it admits a compact generator in the sense of \cref{compact_gen_def}.
\end{proposition}

    For future reference, we give the following results:

\begin{lemma}[{\cite[Lemma II.2.20]{DPS25}}]\label{pseudo_perfect_vs_compact}
	Let $\C \in \Prlo_k$ and let $\Spec (A) \in \dAff_k$.
\begin{enumerate}\itemsep=0.2cm
    \item If $\C$ is smooth, then any pseudo-perfect family in $\C_A$ is compact.
    \item If $\C$ is proper, then any compact object in $\C_A$ is pseudo-perfect.
\end{enumerate}
\end{lemma}

\begin{lemma}[{\cite[Proposition I.4.18]{Por25}}]\label{lemma:tens_prod_finite_type}
    Let $\C, \D \in \Prlo_k$ of finite type.
    Then $\C \otimes_k \D$ is of finite type.
\end{lemma}

\begin{lemma}\label{lemma:pseudo_perf_morph}
    Let $\C \in \Prlo_k$ of finite type and $\Spec(A) \in \dAff_k$.
    Then, the $\infty$-category $\Fun(\Delta^1, \C) \in \Prlo_k$ is of finite type.
    Moreover, an object $(c \to d) \in \Fun(\Delta^1, \C_A) \simeq \Fun(\Delta^1, \C)_A$ is pseudo-perfect if and only if $c,d \in \C_A$ are pseudo-perfect.
\end{lemma}

\begin{proof}
    In the proof of \cite[Lemma 4.4-(2)]{PS23}, it is shown that $\Fun(\Delta^1, \Mod_k)$ is of finite type.
    Hence, the first assertion follows from \cref{lemma:tens_prod_finite_type} and the equivalence $\Fun(\Delta^1, \C) \simeq \Fun(\Delta^1, \Mod_k) \otimes_k \C$.
    For the second assertion, recall that we have an equivalence in $\Mod_A$
\[
\Hom_{\Fun(\Delta^1, \C_A)} (c' \to d', c \to d) \simeq \Hom_{\C_A} (c',c) \times_{\Hom_{\C_A} (c',d)} \Hom_{\C_A} (d',d) \ .
\]
    Recall also the equivalence $\Fun(\Delta^1, \C)^\omega \simeq \Fun(\Delta^1, \C^\omega)$ (see, e.g., the proof of \cite[Lemma 4.4-(1)]{PS23}).
    Then, the if part follows since $\Perf(A)$ is stable under pullback.
    The only if part follows since $0 \xrightarrow{0} e, e \xrightarrow{id} e \in \Fun(\Delta^1, \C)^\omega$ for each $e \in \C^\omega$.
\end{proof}

\subsection{The moduli problem}\label{subsection:moduli_of_objects}
	In this paragraph, we introduce the moduli of objects of an $\infty$-category of \textit{finite type}.
	This construction first appeared in \cite{TV07}, but we will follow the presentation given in \cite[Section 5]{AG14} (see also \cite[Section II.2]{DPS25} and \cite[Section 1.5]{Por25}).
	Notice that our convention for moduli of objects in a category $\C$ differs from the one in \cite{AG14} (but agrees with the one in \cite{TV07, DPS25, Por25}), indeed our moduli of objects of a category $\C$ corresponds to the moduli of objects of the linear dual of $\C$ in \cite{AG14}.
	The two points of view are equivalent and do not affect the results nor the proofs in this paragraph.

\begin{recollection}\label{moduli_of_objects_recollection}
	Let $\C \in \Prlo_k$.
	Consider the presheaf
\begin{align*}
	\widehat{\M}_\C \colon &\dAff_k^{op} \to \mathrm{Spc}\\
	& \Spec(A) \mapsto \C_A^{\simeq}.
\end{align*}
	Since $\C$ is dualizable (\cref{lemma:Prlo_dualizable}), the functor $\C \otimes_k -$ commutes with limits.
	Since $\Mod(-)$ satisfies faithfully flat descent (\cite[Corollary D.6.3.3]{Lur17}), the presheaf $\widehat{\M}_\C $ defines a derived stack over $k$.
	The stack $\widehat{\M}_\C$ is too big and has no chance to be representable.
	The derived moduli of objects of $\C$ is defined as the presheaf 
\begin{align*}
	\M_\C \colon & \dAff_k^{op} \to \mathrm{Spc}\\
	& \Spec(A) \to \C_A^{pp}
\end{align*}
where $\C_A^{pp} \subset \C_A$ is the maximal groupoid spanned by pseudo-perfect families over $A$ in the sense of \cref{pseudo_perfect}.
	Since $\Perf(-)$ satisfies hyper-descent (\cite[Proposition 2.8.4.2-(10)]{Lur18}), the sub-presheaf
\[
\M_\C \subset \widehat{\M}_\C
\]
is a derived substack by \cref{functoriality_tens_prod}.
\end{recollection}

	The main theorem of \cite{TV07} translates as:
\begin{theorem}[{\cite[Theorem 5.8 \& Corollary 5.9]{AG14}}]\label{repr_moduli_of_objects}
	Let $\C \in \Prlo_{k}$ be of finite type.
	Then $\M_\C$ is a locally geometric derived stack locally of finite presentation over $k$.
	Furthermore, the tangent complex at a point $x \colon  \Spec(A) \to \M_\C$ corresponding to a pseudo-perfect family of objects $M_x \in \C_A$ is given by
\[
x^\ast\mathbb{T}_{\M_\C} \simeq \Hom(M_x, M_x)[1] \ .
\]
\end{theorem}

\begin{remark}[{Functoriality of $\M_{(-)}$}]\label{remark:functoriality_moduli_of_objects}
	Let $f \colon \C \to \D$ in $\Prlo_k$ between finite type categories and let $g \colon \D \to \C$ be its right adjoint.
    Assume that $g$ admits a further right adjoint, so that $f$ preserves compact objects.
	Then we have an induced map on moduli of objects
\[
g_{(-)} \colon \M_\D \to \M_\C
\]
	given by $g_A$ on each $\Spec(A) \in \dAff_k$.
	Indeed this map corresponds to the pre-composition with $f$ under the equivalences $\D_A \simeq \Funst_k((\D^{\omega})^{op}, \Mod_A)$ and $\C_A \simeq \Funst_k((\C^{\omega})^{op}, \Mod_A)$.
\end{remark}

\subsection{The moduli of objects of perfect complexes}\label{subsection:moduli_of_objects_Perf}

	We now specialize our discussion in \cref{subsection:moduli_of_objects} to the case of Ind-coherent sheaves on a scheme.
    As proven in Appendix \ref{Appendix:equivalence}, this moduli problem coincides with the stack of perfect complexes with proper support.\\ \indent
	Recall that throughout all this paper, we have fixed an algebraically closed field $k$ of characteristic zero.
	Let $\X \in dSt_k$ and denote by $p\colon X \to \Spec(k)$ the structure morphism.
	Let $\DCoh(\X)$ be the full subcategory of the $\infty$-category $\QCoh(\X)$ spanned by complexes with bounded coherent cohomology sheaves.
	Consider the Ind-completion $\IndCoh(X)$ of $\DCoh(X)$.
	We refer to \cite{GR17} for a thorough study of Ind-coherent sheaves.
	Notice that if $\X = X$ is a smooth scheme, then $\DCoh(X) = \Perf(X)$, hence $\IndCoh(X) = \QCoh(X)$ by Thomason--Trobaugh \cite{TT90}.\\ \indent
	Recall that all the schemes considered in this paper are quasi-compact quasi-separated by assumption.
    The next result of Efimov \cite{Efi20} is crucial for us:
	
\begin{theorem}[{\cite[Theorem 1.4]{Efi20}}]\label{thm:IndCoh_finite}
	Let $X$ be a separated scheme locally of finite presentation over $k$.
	Then $\IndCoh(X) \in \Prlo_k$ is of finite type, hence smooth.
\end{theorem}

    By combining \cref{thm:IndCoh_finite} with \cref{repr_moduli_of_objects}, we get the following result:
	
\begin{corollary}\label{coro:moduli_IndCoh}
	Let $X$ be a separated scheme locally of finite presentation over $k$.
	Then $\M_{\IndCoh(X)}$ is a locally geometric derived stack locally of finite presentation over $k$.
	Furthermore, the tangent complex at a point $x \colon  \Spec(A) \to \M_{\IndCoh(X)}$ corresponding to a pseudo-perfect family of objects $F_x \in \IndCoh(X)_A$ is given by
\[
x^\ast\mathbb{T}_{\M_{\IndCoh(X)}} \simeq \Hom_{\IndCoh(X)_A}(F_x, F_x)[1] \ .
\]
\end{corollary}

\begin{definition}\label{def:moduli_Perf}
	Let $X$ be a scheme over $k$.
	The moduli of perfect complexes on $X$ is
\begin{align*}
\mathbf{Perf}(X) \coloneqq \mathbf{Map}(X, \mathbf{Perf}) \colon & \dAff_k^{op} \to \Spc \\
& \Spec(A) \mapsto  \mathbf{Map}(X_A, \mathbf{Perf}) \simeq \Perf(X_A)^\simeq \ .
\end{align*}
\end{definition}

\begin{definition}\label{def:moduli_Perf_prop}
    In the setting of \cref{def:moduli_Perf}, we let
\[
\mathbf{Perf}_{prop}(X) \subset \mathbf{Perf}(X)
\]
be the substack whose $A$-points are given by the full $\infty$-subgroupoid of $\Perf(X_A)^\simeq$ spanned by those perfect complexes on $X_A$ with support proper over $A$ (\cref{def:proper_support_QCoh}).
\end{definition}

\begin{remark}
    In the setting of \cref{def:moduli_Perf_prop}, if furthermore $X$ is proper, we have an equivalence $\mathbf{Perf}(X) \simeq \mathbf{Perf}_{prop}(X)$.
\end{remark}

    The following result gives a comparison between the stack of perfect complexes with proper support and the moduli of objects in $\IndCoh(X)$.
    
\begin{theorem}[{\cref{thm:comparison_moduli_Perf_DCoh}}]\label{thm:comparison_Perf_corpus}
    Let $X$ be a separated scheme locally of finite presentation over $k$.
    Then, Fourier--Mukai transforms (\cref{construction:Perf_to_Coh}) give an equivalence
\[
\mathbf{Perf}_{prop}(X) \to \M_{\IndCoh(X)} \ .
\]
    In particular, the stack $\mathbf{Perf}_{prop}(X)$ is locally geometric and locally of finite presentation over $k$.
    Moreover, if $X$ is proper, then $\mathbf{Perf}(X) \simeq \M_{\IndCoh(X)}$.
\end{theorem}

\section{Calabi--Yau structures}\label{section:CY_structures}

	In this paragraph, we introduce the notion of relative left Calabi--Yau structures as introduced by Brav--Dyckerhoff \cite{BD19}.
	The main result of this section is \cref{prop:CY_strong}, which implies the existence of Lagrangian structures of \cref{coro:Lagr_strong}.
	The latter is all we shall need in the rest of the paper.
\begin{definition}
	Let $\C \in \mathrm{Pr}^{\mathrm{L}}_k$ dualizable and $(\ev_\C, \coev_\C)$ a duality datum for $\C$.
	Then, the Hochschild homology of $\C$ is defined as the trace of $id_\C$, that is, it is the image of $k$ through the composition
\[
\Mod_k \xrightarrow{\coev_\C} \C \otimes \C^\vee \xrightarrow{\ev_\C} \Mod_k \ .
\]
\end{definition}

\begin{remark}\label{rem:functoriality_HH}
	Let $\C \in \mathrm{Pr}^{\mathrm{L}}_k$ dualizable.
	The formalism of traces of \cite{HSS17} endows $\HH(\C)$ with an action of the circle $S^1$.
	Moreover, Hochschild homology is functorial with respect to morphisms $f\colon \C \to \D$ in $\mathrm{Pr}^{\mathrm{L}}_k$ whose right adjoint admits a further right adjoint.
	In this case, we denote by $\HH(f) \colon \HH(\C) \to \HH(\D)$ the induced morphism, which moreover is $S^1$-equivariant.
	For a thorough account, see \cite[Section 2.5]{Chr23}.
\end{remark}

\begin{definition}\label{def:negative_cyclic_homology}
	Let $\C \in \mathrm{Pr}^{\mathrm{L}}_k$ dualizable.
	The negative cyclic homology of $\C$, denoted by $\HH(\C)^{S^1}$ is the homotopy fixed points of the $S^1$-action on $\HH(\C)$ in \cref{rem:functoriality_HH}.
	It comes equipped with a canonical map in $\Mod_k$
\[
\HH(\C)^{S^1} \to \HH(\C) \ .
\]
\end{definition}

\begin{definition}\label{def:relative_HH}
	Let $\C, \D \in \mathrm{Pr}^{\mathrm{L}}_k$ dualizable and $f \colon \C \to \D$ in $\mathrm{Pr}^{\mathrm{L}}_k$ such that its right adjoint admits a further right adjoint.
\begin{enumerate}\itemsep=0.2cm
	\item The relative Hochschild homology of $f$ is
\[
\HH_f \coloneqq \cof(\HH(\C) \xrightarrow{\HH(f)} \HH(\D)) \in \Mod_k \ .
\]
	When $f$ is clear from the context, we will simply denote $\HH_f$ by $\HH(\D, \C)$.
	\item the relative negative cyclic Hochschild homology of $f$ is 
\[
\HH_f^{S^1} \coloneqq \cof\left(\HH(\C)^{S^1} \to \HH(\D)^{S^1}\right) \in \Mod_k \ .
\]
	When $f$ is clear from the context, we will simply denote $\HH_f^{S^1}$ by $\HH(\D, \C)^{S^1}$.
\end{enumerate}

\end{definition}

\begin{recollection}\label{recollection:inverse_Serre_functor}
	Let $\C \in \Prlo_k$ smooth and $(\ev_\C, \coev_\C)$ be a duality datum for $\C$.
	Then, the functor $\ev_\C$ admits a left adjoint $\ev_\C^{\mathrm{L}} \colon \Mod_k \to \C^\vee \otimes_k \C$, which is automatically $k$-linear.
	Under the equivalence
\[
\Fun_k(\Mod_k, \C^\vee \otimes_k \C) \simeq \Fun_k(\C, \C)
\]
the functor $\ev_\C^{\mathrm{L}}$ corresponds to a $k$-linear endofunctor $Id_\C^! \colon \C \to \C$.
\end{recollection}

\begin{definition}
	Let $\C \in \Prlo_k$ smooth.
	The inverse Serre functor of $\C$ is the $k$-linear endofunctor $Id_\C^! \colon \C \to \C$ of \cref{recollection:inverse_Serre_functor}.
\end{definition}

	The following result relates the inverse Serre functor to Hochschild homology.
	
\begin{lemma}[{\cite[Lemma 2.29]{Chr23}}]\label{lemma:equivalent_def_HH}
	Let $\C \in \Prlo_k$ smooth.
	Then, there exists an equivalence in $\Mod_k$
\[
\HH(\C) \simeq \Fun_k(\Id^!_\C, \Id_\C)
\]
\end{lemma}

\begin{recollection}[{\cite[Section 2.6 \& Section 3.1]{Chr23}}]\label{rem:cyle_to_diagrams_HH}
	Let $\C \in \Prlo_k$ smooth.
	By \cref{lemma:equivalent_def_HH}, a cycle $\sigma \colon k[n] \to \HH(\C)$ corresponds to a $k$-linear morphism
\begin{equation}\label{eq:morph_CY_def}
\widetilde{\sigma} \colon \Id^!_\C[n] \to \Id_\C
\end{equation}
	Moreover, the equivalence of \cref{lemma:equivalent_def_HH} is functorial with respect to left adjoint functors $f \colon \C \to \D$ whose right adjoint $f^r$ admits a further right adjoint.\\ \indent
    Let $f \colon \C \to \D$ in $\mathrm{Pr}^{\mathrm{L}}_k$ between smooth categories and assume that right adjoint $f^r$ to $f$ admits a further right adjoint.
	Consider a relative cycle $\sigma \colon k[n] \to \HH(\D, \C)$, inducing a cycle $\partial\sigma \colon k[n] \to \HH(\C)[1]$.
	By \cref{lemma:equivalent_def_HH}, the relative cycle $\sigma$ corresponds to a diagram of fiber/cofiber sequences in $\Fun_k(\D, \D)$
\begin{equation}\label{eq:triangle_CY_def}
\begin{tikzcd}[sep=small]
	{\Id_\D^!} && {F_!(\Id_\C^!)} && \cof \\
	\\
	\fib && {F_!(\Id_\C)[1-n]} && {\Id_\D[1-n]}
	\arrow["u", from=1-1, to=1-3]
	\arrow["{\widetilde{\sigma}'}"', from=1-1, to=3-1]
	\arrow[from=1-3, to=1-5]
	\arrow["{\widetilde{\partial\sigma}}", from=1-3, to=3-3]
	\arrow["{\widetilde{\sigma}}", from=1-5, to=3-5]
	\arrow[from=3-1, to=3-3]
	\arrow["{cu[1-n]}"', from=3-3, to=3-5]
\end{tikzcd}
\end{equation}
where $F_! \coloneqq f \circ (-) \circ f^r$, $cu$ is the counit of the adjunction $f \dashv f^r$ and $u$ is defined in \cite[Construction 2.31]{Chr23}.
\end{recollection}

	We are now ready to introduce the main notions of this paragraph.
	By \cref{thm:CY_lagr} below, these are non-commutative analogues of shifted symplectic forms and shifted Lagrangian morphisms.
	
\begin{definition}\label{def:CY_structure}
	Let $\C, \D \in \Prlo_k$ smooth and $f \colon \C \to \D$ in $\mathrm{Pr}^{\mathrm{L}}_k$ such that its right adjoint admits a further right adjoint.
\begin{enumerate}\itemsep=0.2cm
	\item We say that a cycle $\sigma \colon k[n] \to \HH(\C)$ is a weak left Calabi--Yau structure of dimension $n$ if the induced morphism $\widetilde{\sigma} \colon \Id^!_\C[n] \to \Id_\C$ of \eqref{eq:morph_CY_def} is an equivalence.
	\item A left Calabi--Yau structure of dimension $n$ is a cycle $\sigma \colon k[n] \to \HH(\C)^{S^1}$ whose induced cycle $k[n] \to \HH(\C)^{S^1} \to \HH(\C)$ is a weak Calabi--Yau structure of dimension $n$.
	\item We say that a relative cycle $\sigma \colon k[n] \to \HH(\D, \C)$ is a weak relative left Calabi--Yau structure of dimension $n$ if $\partial \sigma \colon k[n-1] \to \HH(\C)$ is a weak left Calabi--Yau structure of dimension $n-1$, and the morphism of fiber/cofiber sequences in \eqref{eq:triangle_CY_def} is an equivalence.
	\item A relative left Calabi--Yau structure of dimension $n$ is a relative cycle $\sigma \colon k[n] \to \HH(\D, \C)^{S^1}$ whose induced relative cycle $k[n] \to \HH(\D, \C)^{S^1} \to \HH(\D, \C)$ is a weak relative left Calabi--Yau structure of dimension $n$.
\end{enumerate}
\end{definition}

	For future reference, we give the following example of relative left Calabi--Yau structures, arising from the category of morphisms in a Calabi--Yau category.
	
\begin{proposition}[{\cite[Corollary 6.5]{BD21}}]\label{prop:CY_from_morphisms}
	Let $\C \in \Prlo_k$ smooth equipped with a Calabi--Yau structure $\sigma$ of dimension $d$.
	Then, the functor
\[
((-) \to 0, id_{(-)}, 0 \to (-)) \colon \bigsqcup_{i=1}^{3} \C \to \Fun(\Delta^1, \C)
\]
carries a relative left Calabi--Yau structure of dimension $d+1$ with respect to the Calabi--Yau structure $(\sigma, -\sigma, \sigma)$ on $\bigsqcup_{i=1}^{3} \C$.
\end{proposition}

	The main result of \cite{BD21} relates absolute (respectively, relative) left Calabi--Yau structures to shifted symplectic (respectively, Lagrangian) structures.
	
\begin{theorem}[{\cite[Theorem 5.6]{BD21}}]\label{thm:CY_lagr}\hfill
\begin{enumerate}[(1)]\itemsep=0.2cm
	\item Let $\C \in \Prlo$ smooth.
	Then, a left $d$-Calabi--Yau structure on $\C$ induces a $(2-d)$ shifted symplectic structure on $\M_\C$. 
	\item Let $\C, \D \in \Prlo_k$ smooth and $f \colon \C \to \D$ in $\mathrm{Pr}^{\mathrm{L}}_k$ such that its right adjoint $f^r$ admits a further right adjoint.
	Then, a relative left $d$-Calabi--Yau structure on $f$ induces a $(2-d)$-shifted Lagrangian structure on $f^r \colon \M_\D \to \M_\C$.
\end{enumerate}
\end{theorem}

	Combining the relative left Calabi--Yau structures from \cref{prop:CY_from_morphisms} with \cref{thm:CY_lagr}, we have:
	
\begin{corollary}\label{coro:lagr_from_morphims}
	Let $\C \in \Prlo_k$ equipped with a Calabi--Yau structure $\sigma$ of dimension $d$.
	Then $\M_\C$ carries a $(2-d)$ shifted symplectic structure $\omega$ and the morphism in $\dSt_k$
\[
(\fib, \ev_0, \ev_1) \colon \M_{\Fun(\Delta^1, \C)} \to \prod_{i=1}^{3} \M_\C
\]
carries a $(1-d)$-shifted Lagrangian structure with respect to the $(2-d)$-shifted symplectic structure $(\omega, -\omega, \omega)$ on $\prod_{i=1}^{3} \M_\C$, where $\ev_0$ (respectively, $\ev_1$, respectively $\fib$) sends a morphism between pseudo-perfect objects to its source (respectively, target, fiber).
\end{corollary}

\subsection{Algebraic geometry}\label{subsection:CY_algebraic_geometry}

	We now specialize our discussion to Calabi--Yau structures arising from algebraic geometry.
	The Calabi--Yau structures considered in this paragraph are a reflection of Grothendieck--Serre duality.
	We will start by recalling the results of \cite[Section 5.2]{BD19}.\\ \indent
	By \cref{thm:IndCoh_finite}, it makes sense to study Calabi--Yau structures on $\IndCoh(X)$ for $X$ a separated scheme locally of finite presentation over $k$.
	We will see that $\IndCoh(X)$ carries a left Calabi--Yau structure if and only if $X$ is Calabi--Yau, i.e., if $X$ has trivial canonical bundle.\\ \indent
	Let $f \colon X \to Y$ be a proper morphism between separated schemes locally of finite presentation over $k$.
	Then, we have an adjunction
\[
f_\ast \colon \IndCoh(X) \leftrightarrows \IndCoh(Y) \colon f^! \ .
\]
	Since $f$ is proper, we have that $f_\ast$ preserves compact objects, that is, we have an inclusion $f_\ast(\DCoh(X)) \subset \DCoh(Y)$, so that the right adjoint $f^!$ admits a further right adjoint.
	Notice that since $X$ is separated, the above can be applied to the diagonal $\Delta \colon X \to X \times_k X$, hence giving an adjunction $\Delta_\ast \dashv \Delta^!$.

\begin{recollection}\label{recollection:upper_shriek_compactification}
    The functor $f^! \colon \IndCoh(Y) \to \IndCoh(X)$ is more generally defined for arbitrary morphisms between (derived) schemes locally (almost) of finite presentation in \cite{GR17}.
    The following special case will be enough for us.
    Given a separated scheme locally of finite presentation over $k$, choose a compactification $\overline{X}$ of $X$ over $k$
\[
\begin{tikzcd}[sep=small]
	X && {\overline{X}} \\
	\\
	& {\Spec(k)}
	\arrow["j", hook, from=1-1, to=1-3]
	\arrow["p"', from=1-1, to=3-2]
	\arrow["q", from=1-3, to=3-2]
\end{tikzcd}
\]
where $j$ is an open immersion and $q$ is proper.
    By \cite[Chapter 5, § 4.3]{GR17}, we have
\[
p^! \simeq j^\ast \circ q^! \colon \Mod_k \to \IndCoh(X) \ .
\]
\end{recollection}
	
\begin{definition}\label{def:dualizing_complex}
	Let $X$ be a separated scheme locally of finite presentation over $k$ and let $p \colon X \to \Spec(k)$ be the structure morphism.
	The dualizing complex of $X$ is $\omega_X^\bullet \coloneqq p^!(k)$.
\end{definition}

\begin{remark}
	Notice that in \cref{def:dualizing_complex}, keeping the notation of \cref{recollection:upper_shriek_compactification}, we could have also considered the functor $j^\ast \circ q^{!, \QCoh} \colon \Mod_k \to \QCoh(X)$, where $q^{!, \QCoh}$ is the right adjoint of $q_\ast \colon QCoh(\overline{X}) \to \Mod_k$.
	By \cite[Chapter 4, Lemma 5.1.9]{GR17}, this would not change the definition of $\omega_X^\bullet$ (at least after choosing a compactification of $X$ as in \cref{recollection:upper_shriek_compactification}).
	Nevertheless, for the functoriality of Hochschild homology along pushforward, we need to ensure that its right adjoint admits a further right adjoint.
	This holds for pushforward along proper morphisms in $\IndCoh$, but it does not hold in $\QCoh$ in general.
	This is the main advantage of working with Ind-coherent sheaves instead of quasi-coherent sheaves.
\end{remark}

\begin{recollection}\label{recollection:auto_duality_IndCoh}
	Ind-coherent sheaves satisfy a categorical Künneth formula, that is, we have an equivalence
\begin{equation}\label{eq:Kunneth}
\IndCoh(X) \otimes_k \IndCoh(Y) \xrightarrow{\sim}\IndCoh(X \times_k Y)
\end{equation}
which is given by the Ind-completion of the functor sending $F \in \DCoh(X), G \in \DCoh(Y)$ to $\pi_X^\ast F \otimes_{\OO_{X \times_k Y}} \pi_Y^\ast G$, where $\pi_X \colon X \times_k Y \to X$ is the projection, and similarly for $\pi_Y$.
	Using the equivalence in \eqref{eq:Kunneth}, Gaitsgory--Rozenblyum \cite[Chapter 5, § 4.3]{GR17} show that $\IndCoh(X)$ is auto-dual with duality datum given by
\begin{align}\label{eq:duality_IndCoh}
\begin{split}
\ev_X \colon \IndCoh(X \times_k X) \xrightarrow{\Delta^!} \IndCoh(X) \xrightarrow{p_\ast} \Mod_k \ , \\
\coev_X \colon  \Mod_k\xrightarrow{p^!} \IndCoh(X) \xrightarrow{\Delta_\ast} \IndCoh(X \times_k X) \ .
\end{split}
\end{align}
\end{recollection}

	The duality datum in \eqref{eq:duality_IndCoh} gives the following formula for Hochschild homology:
	
\begin{lemma}[{\cite[Lemma 5.8]{BD19}}]\label{lemma:HH_IndCoh}
	For $X$ a separated scheme locally of finite presentation over $k$, there is an equivalence
\[
\HH(X) \simeq \Gamma(X, \Delta^!\Delta_\ast \omega_X^\bullet) \simeq \Hom_{X\times_k X}(\Delta_\ast \OO_X, \Delta_\ast \omega_X^\bullet) \ .
\]
\end{lemma}

	Furthermore, Brav--Dyckerhoff give also an explicit description for the functoriality of Hochschild homology along proper morphisms:

\begin{lemma}[{\cite[Lemma 5.9]{BD19}}]\label{lemma:functoriality_HH_proper}
	Let $f \colon X \to Y$ be a proper morphism between separated schemes locally of finite presentation over $k$.
	Under the equivalence of \cref{lemma:HH_IndCoh}, the natural map $\HH(f_\ast) \colon \HH(\IndCoh(X)) \to \HH(\IndCoh(Y))$ corresponds to the composition
\[
\Hom_{X\times_k X}(\Delta_\ast \OO_X, \Delta_\ast \omega_X^\bullet) \to \Hom_{Y\times_k Y}(\Delta_\ast f_\ast \OO_X, \Delta_\ast f_\ast \omega_X^\bullet) \to \Hom_{Y\times_k Y}(\Delta_\ast \OO_Y, \Delta_\ast \omega_Y^\bullet) \ ,
\]
where the first arrow is given by pushforward along $(f, f) \colon X \times_k X \to Y \times_k Y$ and the second map is given by pre-composition with $\OO_Y \to f_\ast \OO_X$ and post-composition with $f_\ast \omega_X^\bullet \simeq f_\ast f^! \omega_Y^\bullet \to \omega_Y^\bullet$.
\end{lemma}

	Our objects of study are allowed to be singular.
	Still, we will impose that the singularities are mild.

\begin{definition}\label{def:Goresentein_CY}
	Let $X$ be a separated scheme locally of finite presentation over $k$.
\begin{enumerate}\itemsep=0.2cm
	\item We say that $X$ is a Gorenstein scheme of dimension $d$ if $\omega_X^\bullet[-d]$ is a line bundle.
	\item We say that $X$ is Calabi--Yau of dimension $d$ if it is Gorenstein of dimension $d$ and there exists a trivialization $\omega_X^\bullet[-d] \simeq \OO_X$.
\end{enumerate}
\end{definition}


	\cref{def:Goresentein_CY}-(2) is justified by the following result:
	
\begin{theorem}[{\cite[Proposition 5.12]{BD19}}]\label{thm:CY_IndCoh}
	Let $X$ be a Gorenstein scheme of dimension $d$.
	Then giving a left Calabi--Yau structure on $\IndCoh(X)$ is equivalent to giving a trivialization $\OO_X \simeq \omega_X^\bullet [-d]$.
\end{theorem}

    By \cref{thm:comparison_Perf_corpus}, \cref{thm:CY_lagr}, and \cref{thm:CY_IndCoh}, we get:

\begin{corollary}\label{coro:sympl_Calbi_Yau}
    Let $X$ be a Calabi--Yau scheme of dimension $d$.
    Then $\mathbf{Perf}_{prop}(X)$ carries a $(2-d)$-shifted symplectic structure.
\end{corollary}

\begin{definition}\label{def:strong_CY_scheme}
	We say that a Gorenstein scheme $X$ of dimension $d$ is \textit{strong} Calabi--Yau if $X$ is proper, connected, Calabi--Yau of dimension $d$, and $\mathrm{H}^i(X, \OO_X)$ vanishes for $i \neq 0, d$ and has dimension $1$ for $i=0, d$.
\end{definition}

\begin{remark}\label{rem:lichtenbaum_properness}
    Let $X$ be an irreducible,  Calabi--Yau scheme of dimension $d$, and assume that $\mathrm{H}^i(X, \OO_X)$ vanishes for $i \neq 0, d$ and has dimension $1$ for $i=0, d$.
    Then $X$ is automatically strong Calabi--Yau.
    Indeed, connectedness follows from $X$ being irreducible, and properness follows by the Grothendieck--Lichtenbaum theorem \cite[\href{https://stacks.math.columbia.edu/tag/0G5F}{Theorem 0G5F}]{stacks-project}.
\end{remark}

\begin{remark}\label{rem:CY_scheme_proepr}
    In the setting of \cref{def:strong_CY_scheme}, if furthermore $X$ is reduced, then by Grothendieck--Serre duality and the triviality of $\omega_X^\bullet[-d]$, we have that $\mathrm{H}^d(X, \OO_X) = \mathrm{H}^0(X, \OO_X) = k$.
    Hence, the last condition in \cref{def:strong_CY_scheme} becomes redundant in this case.
\end{remark}

\begin{remark}
	In \cref{def:strong_CY_scheme}, if $X$ is smooth of dimension $2$, we recover the definition of a K3 surface.
\end{remark}

	We are now ready to prove the main result of this section:

\begin{proposition}\label{prop:CY_strong}
	Let $X$ be a proper Gorenstein scheme of dimension $d > 0$ and let $p \colon X \to \Spec(k)$ be the structure morphism.
	Then, the functor
\[
p_\ast \colon \IndCoh(X) \to \Mod_k
\]
carries a relative left Calabi--Yau structure of dimension $d+1$ if and only if $X$ is a strong Calabi--Yau scheme.
\end{proposition}

\begin{proof}
	By \cite[Lemma 5.10]{BD19}, the canonical morphism $\HH_{d+1}(\Mod_k, \IndCoh(X))^{S^1} \to \HH_{d+1}(Mod_k, \IndCoh(X))$ is an equivalence.
	Hence it is enough to show that the functor $p_\ast \colon \IndCoh(X) \to \Mod_k$ admits a \textit{weak} relative left Calabi--Yau structure of dimension $d+1$ if and only if $X$ is a strong Calabi--Yau scheme of dimension $d$.
	By \cref{rem:cyle_to_diagrams_HH} and \cref{lemma:functoriality_HH_proper}, and since the diagonal of $\Spec(k)$ is the identity $id_{\Spec(k)}$, giving a weak relative left Calabi--Yau structure of dimension $d+1$ on $p_\ast$ is equivalent to giving:
\begin{enumerate}\itemsep=0.2cm
	\item an equivalence $\varphi \colon \OO_X \stackrel{\sim}{\to} \omega_X^\bullet[-d]$;
	\item a homotopy to zero of the composition
\begin{equation}\label{eq:diagram_prop_CY}
\begin{tikzcd}[sep=small]
	k && {p_\ast \OO_X } && \cof \\
	\\
	\fib && {p_\ast \omega_X^\bullet[-d]} && {k[-d]}
	\arrow[from=1-1, to=1-3]
	\arrow[dashed, from=1-1, to=3-1]
	\arrow[from=1-3, to=1-5]
	\arrow[from=1-3, to=3-3]
	\arrow["{p_\ast(\varphi)}"{description}, shift left=5, draw=none, from=1-3, to=3-3]
	\arrow[dashed, from=1-5, to=3-5]
	\arrow[from=3-1, to=3-3]
	\arrow[from=3-3, to=3-5]
\end{tikzcd}
\end{equation}
such that the induced dotted arrows are equivalences.
\end{enumerate}
	By \cref{thm:CY_IndCoh}, item $(1)$ is equivalent to $X$ being Calabi--Yau of dimension $d$.
	Hence, we are reduced to showing that, assuming item $(1)$, a homotopy as in item $(2)$ exists if and only if $X$ is a strong Calabi--Yau scheme.\\ \indent
	If 
\[
k \not\simeq \mathrm{H}^0\Gamma(X, \OO_X) \simeq \tau_{\geq 0}\Gamma(X, \OO_X) \ ,
\]
then the right vertical arrow in \eqref{eq:diagram_prop_CY} cannot be an equivalence.
	Hence, we can assume $\dim_k H^0(X, \OO_X) = 1$, in which case the top fiber/cofiber sequence in \eqref{eq:diagram_prop_CY} is canonically identified with the fiber/cofiber sequence
\[
\tau_{\geq 0}\Gamma(X, \OO_X) \to \Gamma(X, \OO_X) \to \tau_{\leq - 1}\Gamma(X, \OO_X) \ .
\]
	Similarly, the composition of the central vertical map and bottom right map in \eqref{eq:diagram_prop_CY} is identified with the canonical morphism $p_\ast \OO_X \to \mathrm{H}^d \Gamma(X, \OO_X) \simeq \tau_{\leq -d} \Gamma(X, \OO_X)$.
	Hence, the composition \eqref{eq:diagram_prop_CY} is identified with the composition
\[
\tau_{\geq 0}\Gamma(X, \OO_X) \to \Gamma(X, \OO_X) \to \tau_{\leq -d} \Gamma(X, \OO_X) \ ,
\]
which manifestly has a homotopy to zero, as the second arrow factors through the map $\Gamma(X, \OO_X) \to \tau_{\leq -1} \Gamma(X, \OO_X)$.
	Under the above identifications, \eqref{eq:diagram_prop_CY} can be rewritten as
\[
\begin{tikzcd}[sep=small]
	k && {p_\ast \OO_X } && {\tau_{\leq -1} \Gamma(X, \OO_X)} \\
	\\
	\fib && {p_\ast \omega_X^\bullet[-d]} && {\tau_{\leq -d} \Gamma(X, \OO_X)}
	\arrow[from=1-1, to=1-3]
	\arrow[dashed, from=1-1, to=3-1]
	\arrow[from=1-3, to=1-5]
	\arrow[from=1-3, to=3-3]
	\arrow["\sim"{marking, allow upside down}, shift left, draw=none, from=1-3, to=3-3]
	\arrow["\psi", dashed, from=1-5, to=3-5]
	\arrow[from=3-1, to=3-3]
	\arrow[from=3-3, to=3-5]
\end{tikzcd}
\]
where $\psi \colon \tau_{\leq -1} \Gamma(X, \OO_X) \to \tau_{\leq -d} \Gamma(X, \OO_X)$ is the canonical morphism.
	By definition, the morphism $\psi$ is an equivalence if and only if $X$ is \textit{strong} Calabi--Yau.
\end{proof}

	Combining \cref{thm:CY_lagr} with \cref{remark:functoriality_moduli_of_objects} and \cref{thm:comparison_Perf_corpus}, we have:
	
\begin{corollary}\label{coro:Lagr_strong_exceptional}
	Let $X$ be a strong Calabi--Yau scheme of dimension $d > 0$.
	Then, the morphism in $\dSt_k$
\[
p^! \colon \mathbf{Perf}(k) \to \mathbf{Perf}(X) 
\]
carries a $(1-d)$-shifted Lagrangian structure.
\end{corollary}

\begin{remark}
    \cref{prop:CY_strong} has two main advantages with respect to the commutative statement of \cref{coro:Lagr_strong_exceptional} (and \cref{coro:Lagr_strong} below).
	Firstly, it is a finer statement.
	Secondly, invoking \cite[Lemma 5.10]{BD19}, in the proof we could reduce to showing the existence of a weak relative left Calabi--Yau structure, hence avoiding the need of taking care of the circle action on Hochschild homology.
	The second point translates, in the commutative setting, into the fact that it is enough to construct a non-degenerate isotropic structure for the underlying $2$-form of the symplectic form on $\mathbf{Perf}(X)$.
    One can also give a direct, commutative proof of \cref{coro:Lagr_strong}, see Appendix \ref{appendix:the_case_of_Perf}.
\end{remark}

\begin{corollary}\label{coro:Lagr_strong}
    Let $X$ be a strong Calabi--Yau scheme of dimension $d > 0$.
	Then, the morphism in $\dSt_k$
\[
p^\ast \colon \mathbf{Perf}(k) \to \mathbf{Perf}(X) 
\]
carries a $(1-d)$-shifted Lagrangian structure.
\end{corollary}

\begin{proof}
    Since $\omega_X^\bullet[-d] \simeq \OO_X$ by assumption, we have an equivalence of functors
\[
p^! \simeq p^\ast [d] \colon \Perf(k) \to \Perf(X) \ .
\]
    More generally, we have
\[
p^! \otimes_k \Id_{Mod_A} \simeq p^\ast [d] \otimes_k \Id_{Mod_A} \simeq p_A^\ast [d] \colon \Perf(A) \to \Perf(X_A) \ ,
\]
where $p_A \colon X_A \to \Spec(A)$ is the projection.
    Hence, the result follows from \cref{coro:Lagr_strong_exceptional} as the claim is stable under shifts.
\end{proof}

	For further reference, we also recall examples of relative left Calabi--Yau and Lagrangian morphisms arising from anti-canonical divisors \cite{BD19, BD21}.
	
\begin{theorem}[{\cite[Theorem 5.13]{BD19}}]\label{thm:CY_divisors}
	Let $Y$ be a Gorenstein scheme over $k$ of dimension $d+1$ and fix a global section $s \in \mathrm{H}^0\Gamma(Y, \omega_Y^{-1})$, where $\omega_Y \coloneqq \omega_Y^\bullet[-d-1]$.
	Assume that the zero locus $i \colon X \hookrightarrow Y$ of $s$ has dimension $d$.
    Then $X$ is Calabi--Yau of dimension $d$ and the functor
\[
i_\ast \colon \IndCoh(X) \to \IndCoh(Y)
\]
carries a relative left Calabi--Yau structure of dimension $d+1$.
\end{theorem}

	As an immediate consequence of \cref{thm:comparison_Perf_corpus}, \cref{thm:CY_lagr} and \cref{thm:CY_divisors}, we have:
	
\begin{corollary}[{\cite[Corollary 6.2]{BD21}}]\label{thm:Lagr_divisors}
	In the setting of \cref{thm:CY_divisors}, the derived stack $\mathbf{Perf}_{prop}(X)$ carries a $(2-d)$-shifted symplectic structure and the morphism in $\dSt_k$
\[
i^! \colon \mathbf{Perf}_{prop}(Y) \to \mathbf{Perf}_{prop}(X)
\]
carries a $(1-d)$-shifted Lagrangian structure.
\end{corollary}

\section{The main result}\label{section:main_results}

\subsection{Composition of Lagrangian cospans}

	For future reference, we recall here a result on composition of Lagrangian cospans of Calaque \cite{Cal15}, and its non-commutative analogue of Brav-Dyckerhoff \cite{BD19}.

\begin{theorem}[{\cite[Theorem 4.4]{Cal15}, \cite[Theorem 6.2]{BD19}}]\label{thm:composition_lagr_morph}
	Consider a diagram in $\dSt_k$
\[
\begin{tikzcd}[sep=small]
	{\mathcal{L}_1 \times_{\mathcal{Y}} \mathcal{L}_2} && {\mathcal{L}_2} && \mathcal{Z} \\
	\\
	{\mathcal{L}_1} && \mathcal{Y} \\
	\\
	{\mathcal{X}}
	\arrow["{p_2}", from=1-1, to=1-3]
	\arrow["{p_1}"', from=1-1, to=3-1]
	\arrow["{f_2}", from=1-3, to=1-5]
	\arrow["{g_2}", from=1-3, to=3-3]
	\arrow["{g_1}"', from=3-1, to=3-3]
	\arrow["{f_1}"', from=3-1, to=5-1]
\end{tikzcd}
\]
	and assume that:
\begin{enumerate}\itemsep=0.2cm
	\item $\mathcal{X}, \mathcal{Y}, \mathcal{Z}$ carry $n$-shifted symplectic forms $\omega_X, \omega_{\mathcal{Y}}, \omega_{\mathcal{Z}}$;
	\item the morphism $(f_1, g_1) \colon \mathcal{L}_1 \to \mathcal{X} \times_k \mathcal{Y}$ carries a $(n-1)$-shifted Lagrangian structure with respect to $(\omega_{\mathcal{X}}, -\omega_{\mathcal{Y}})$, and the morphism $(g_2, f_2) \colon \mathcal{L}_2 \to \mathcal{Y} \times_k \mathcal{Z}$ carries a $(n-1)$-shifted Lagrangian structure with respect to $(\omega_{\mathcal{Y}}, -\omega_{\mathcal{Z}})$.
\end{enumerate}
	Then, the morphism
\[
(f_1 \circ p_1, f_2 \circ p_2) \colon \mathcal{L}_1 \times_{\mathcal{Y}} \mathcal{L}_2 \to \mathcal{X} \times_k \mathcal{Z}
\]
carries a $(n-1)$-shifted Lagrangian structure with respect to $(\omega_{\mathcal{X}}, -\omega_{\mathcal{Z}})$.\\ \indent
	Moreover, if every Lagrangian structure in the assumptions arises via \cref{thm:CY_lagr}, then the Lagrangian structure on $(f_1 \circ p_1, f_2 \circ p_2)$ also arises via \cref{thm:CY_lagr}.
\end{theorem}

\begin{proof}
	The first part of the claim is exactly \cite[Theorem 4.4]{Cal15}.
	The moreover part is \cite[Theorem 6.2]{BD19}.
	See also \cite[Theorem 3.15]{Chr23}.
\end{proof}

\subsection{Shifted Lagrangian morphisms from syzygies}\label{subsection:Lagrangian_IndCoh}

\begin{definition}
    Let $X$ be a separated scheme locally of finite presentation over $k$.
    We define the derived stack over $k$
\begin{align*}
\mathbf{Perf}_{prop}(X)^{\Delta^1} \coloneqq \mathbf{Map}(\Delta^1, \mathbf{Perf}_{prop}(X))\colon & \dAff_k^{op} \to \Spc \\
& \Spec(A) \to \Map(\Delta^1, \Perf_{prop / A}(X_A)) \ .
\end{align*}
\end{definition}

\begin{proposition}\label{prop:Perf_ext}
    Let $X$ be a separated scheme locally of finite presentation over $k$.
    Then $\mathbf{Perf}_{prop}(X)^{\Delta^1}$ is locally geometric and locally of finite presentation over $k$.
\end{proposition}

\begin{proof}
    It follows by \cref{lemma:pseudo_perf_morph}, \cref{repr_moduli_of_objects} and \cref{thm:comparison_Perf_corpus} that we have an equivalence
\[
\mathbf{Perf}_{prop}(X)^{\Delta^1} \simeq \M_{\Fun(\Delta^1, \IndCoh(X))} \ ,
\]
and $\M_{\Fun(\Delta^1, \IndCoh(X))}$ is locally geometric and locally of finite presentation over $k$.
\end{proof}
	
\begin{definition}\label{def:Syz_L}
    Let $X$ be a proper scheme locally of finite presentation over $k$ and $i \colon \mathbf{\mathcal{L}} \to \mathbf{Perf}(X)$ be a morphism of stacks over $k$.
	The derived stack $\widetilde{\mathbf{Syz}}^{\mathcal{L}}(X)$ is defined by the pullback square in $\dSt_k$
\[
\begin{tikzcd}[sep=small]
	{\widetilde{\mathbf{Syz}}^{\mathcal{L}}(X)} && {\mathbf{Perf}(k) \times \mathbf{\mathcal{L}}} \\
	\\
	{\mathbf{Perf}(X)^{\Delta^1}} && {\mathbf{Perf}(X) \times \mathbf{Perf}(X)}
	\arrow["{(\widetilde{\ev}_0, \widetilde{\ev}_1)}", from=1-1, to=1-3]
	\arrow["q"', from=1-1, to=3-1]
	\arrow["{(p^\ast, i)}", from=1-3, to=3-3]
	\arrow["{(\ev_0, \ev_1)}"', from=3-1, to=3-3]
\end{tikzcd}
\]
\end{definition}

	We will need the following result, obtained by specializing \cref{coro:lagr_from_morphims} to $\C = \IndCoh(X)$ and using the equivalences of \cref{thm:comparison_Perf_corpus} and \cref{prop:Perf_ext}:
	
\begin{proposition}\label{prop:Lagrangian_ev_fib}
	Let $X$ be a Calabi--Yau scheme of dimension $d$.
	Then, the morphism in $\dSt_k$
\[
(\ev_0, \ev_1, \fib) \colon \mathbf{Perf}_{prop}(X)^{\Delta^1}\to \mathbf{Perf}_{prop}(X)^{\times 3}
\]
carries a $(1-d)$-shifted Lagrangian structure.
\end{proposition}

	Putting all together, we obtain the main result of the paper:
	
\begin{theorem}\label{thm:main_DCoh}
	Let $X$ be a strong Calabi--Yau scheme of dimension $d > 0$ and $i \colon \mathbf{\mathcal{L}} \to \mathbf{Perf}(X)$ be a $(1-d)$-shifted Lagrangian morphism.
	Then, the morphism of derived stacks
\[
\fib \circ q \colon \widetilde{\mathbf{Syz}}^{\mathbf{\mathcal{L}}}(X) \to \mathbf{Perf}(X)
\]
carries a $(1-d)$-shifted Lagrangian structure.
	Moreover, if the $(1-d)$-shifted Lagrangian structure on $i \colon \mathbf{\mathcal{L}} \to \mathbf{Perf}(X)$ arises from a relative left Calabi--Yau structure of dimension $d+1$ via \cref{thm:CY_lagr}, then the $(1-d)$-shifted Lagrangian structure on $\fib \circ q$ also arises via \cref{thm:CY_lagr}.
\end{theorem}

\begin{proof}
    Consider the diagram in $\dSt_k$ with pullback square
\begin{equation}\label{eq:def_Syz_tilde}
\begin{tikzcd}[sep=small]
	{\widetilde{\mathbf{Syz}}^{\mathcal{L}}(X)} && {\mathbf{Perf}(k) \times \mathbf{\mathcal{L}}} \\
	\\
	{\mathbf{Perf}(X)^{\Delta^1}} && {\mathbf{Perf}(X) \times \mathbf{Perf}(X)} \\
	\\
	{\mathbf{Perf}(X)}
	\arrow["{(\widetilde{\ev}_0, \widetilde{\ev}_1)}", from=1-1, to=1-3]
	\arrow["q"', from=1-1, to=3-1]
	\arrow["{(p^\ast, i)}", from=1-3, to=3-3]
	\arrow["{(\ev_0, \ev_1)}"', from=3-1, to=3-3]
	\arrow["\fib"', from=3-1, to=5-1]
\end{tikzcd}
\end{equation}
	The result then follows directly by \cref{thm:comparison_Perf_corpus}, \cref{coro:Lagr_strong}, \cref{thm:composition_lagr_morph} and \cref{prop:Lagrangian_ev_fib}.
\end{proof}

\begin{example}\label{example:Strong_CY_from_anticanonical_div}
    In the setting of \cref{thm:CY_divisors}, we have a $(1-d)$-shifted Lagrangian morphism $i^! \colon \mathbf{Perf}_{prop}(Y) \to \mathbf{Perf}_{prop}(X)$.
	If we further assume that $d>0$ and $X$ is strong Calabi--Yau (e.g., if $Y$ is a Fano variety by \cite[Theorem 0.1]{KMM92} and \cite[Corollary 2.6]{Voi16}), then the assumptions of \cref{thm:main_DCoh} are satisfied and the morphism
\[
\fib \circ q \colon \widetilde{\mathbf{Syz}}^{\mathbf{Perf}_{prop}(Y)}(X) \to \mathbf{Perf}(X)
\]
is $(1-d)$-shifted Lagrangian.
\end{example}

\section{Comparison with the work of Fantechi--Miró-Roig}\label{section:comparison_with_FMR}

	In this Section, we provide a comparison with the work of Fantechi--Miró-Roig \cite{FMR25, FMR25b}.
    In particular, we show how to recover and generalize their main results on symplectic and Lagrangian morphisms in \textit{loc.cit.} from the results obtained in \cref{section:main_results}.

\subsection{The open substack of simple perfect complexes}

	We introduce the open substacks of simple perfect complexes $\mathbf{Perf}^{sim}(X)$ and simple perfect sheaves $\Spl(X)$.
	The algebraic space of simple perfect sheaves obtained by rigidifying $\Spl(X)$ is the main object of study of \cite{FMR25,FMR25b}.\\ \indent
    
\begin{definition}[{\cite[Remark 3.23]{TV07}, \cite[Definition 5.1]{STV_determinant_perfect_complexes}}]\label{def:Perf_simple}
    Let $X$ be a proper scheme locally of finite presentation over $k$.
    The stack of simple perfect complexes on $X$ is the derived open substack
\[
\mathbf{Perf}^{\simple}(X) \subset \mathbf{Perf}(X)
\]
corresponding to the classical open substack
\[
t_0\mathbf{Perf}^{\simple}(X) \subset t_0\mathbf{Perf}(X)
\]
sending $\Spec(A) \in \Aff_k$ to the full $1$-subgroupoid of $\mathrm{Perf}(X_A)^\simeq$ spanned by perfect complexes $\E \in \Perf(X_A)$ such that, for each $\Spec(A') \to \Spec(A)$ in $\Aff_k$, we have $\Ext^i(\E \otimes_A A', \E \otimes_A A') = 0$ for $i < 0$, and the canonical map $A' \to \Ext^0(\E \otimes_A A', \E \otimes_A A')$ is an isomorphism.
\end{definition}

Thanks to the work of Montagnani \cite{Mon26}, it is possible to explicitly describe the points of the derived stack $\mathbf{Perf}^{\simple}(X)$ of \cref{def:Perf_simple}:
\begin{remark}
    Keeping the same notations as in \cref{def:Perf_simple}, \cite[Section 2]{Mon26} shows that $\mathbf{Perf}^{\simple}(X)(\Spec(A))$ is the full $\infty$-subgroupoid of $\mathrm{Perf}(X_A)^\simeq$ spanned by those perfect complexes $\E \in \Perf(X_A)$ such that the cofiber of the natural map in $\Mod_A$ 
\[
A \to \Hom(\E, \E)
\]
has Tor-amplitude strictly smaller than $0$.
    That the above definition recovers \cref{def:Perf_simple} is the content of \cite[Lemma 2.11 \& Corollary 2.14]{Mon26}.
    We warn the reader that if one attempts to reproduce naively \cref{def:Perf_simple} but with $A$ derived, that is, if one requires the natural map $A \to \tau_{\geq 0}\Hom(\E, \E)$ to be an equivalence, then this condition is not even stable by base change.
    See \cite[Example 2.15]{Mon26}.
\end{remark}

\begin{remark}
    That the substack $\mathbf{Perf}^{\simple}(X) \subset \mathbf{Perf}(X)$ of \cref{def:Perf_simple} is an open substack follows by semicontinuity of dimension.
    Moreover, it is a $1$-Artin stack.
    See \cite[Section 3.4]{TV07}.
\end{remark}

\begin{definition}\label{def:Perf_Tor_amplitude}
    Let $X$ be a proper scheme locally of finite presentation over $k$.
    The stack of perfect complexes of Tor-amplitude $[0,0]$ on $X$ is the open substack
\[
\mathbf{Perf}^{[0,0]}(X) \subset \mathbf{Perf}(X)
\]
sending $\Spec(A)$ to the maximal $\infty$-subgroupoid of $\Perf(X_A)^\simeq$ spanned by those $\E \in \Perf(X_A)$ of Tor-amplitude $[0,0]$, that is, for which $\E \otimes_{\OO_{X_A}} p_A^\ast(M) \in \QCoh(X_A)^\heartsuit$ for each $M \in \Mod_A^\heartsuit$, where $p_{A} \colon X_A \to \Spec(A)$ is the projection morphism.
\end{definition}

\begin{remark}\label{rem:Perf00_1_Artin}
    That the substack $\mathbf{Perf}^{[0,0]}(X) \subset \mathbf{Perf}(X)$ of \cref{def:Perf_Tor_amplitude} is an open substack follows by \cite[Proposition 6.1.4.5]{Lur18} and the fact that $X$ is quasi-compact.
    Moreover, it is a $1$-Artin stack by \cite[Lemma 2.19]{TV07}.
\end{remark}

\begin{definition}\label{def:simple_sheaves}
    Let $X$ be a proper scheme locally of finite presentation over $k$.
    The stack of simple sheaves on $X$ is the open substack
\[
\mathbf{Spl}(X) \coloneqq\mathbf{Perf}^{\simple}(X) \cap \mathbf{Perf}^{[0,0]}(X) \subset \mathbf{Perf}(X) \ .
\]
\end{definition}

\begin{remark}
    Notice that the $A$-points of the stack $\mathbf{Spl}(X)$ do not need to be \textit{sheaves}, that is, complexes that are concentrated in degree $0$.
    Nevertheless, this holds whenever $\Spec(A)$ is a classical affine scheme, so that the points of the truncation $t_0(\Spl(X))$ are indeed perfect sheaves.
\end{remark}

	We now give two equivalent definitions of the derived stack $\widetilde{\mathbf{Syz}}^{\mathcal{L}}(X)$ of \cref{def:Syz_L}.
\begin{lemma}\label{lemma:equiv_def_Syz_L}
	Let $X$ be a proper scheme locally of finite presentation over $k$.
    The derived stack $\widetilde{\mathbf{Syz}}^{\mathcal{L}}(X)$ is equivalent to the one defined by the following diagrams in $\dSt_k$ with pullback squares:
\[
\begin{tikzcd}[sep=small]
	{\widetilde{\mathbf{Syz}}^{\mathcal{L}}(X)} && {\mathbf{Perf}(X)^{\Delta^1, \mathbf{\mathcal{L}}}} && {\mathbf{\mathcal{L}}} \\
	\\
	{\widetilde{\mathbf{Syz}}(X)} && {\mathbf{Perf}(X)^{\Delta^1}} && {\mathbf{Perf}(X)} \\
	\\
	{\mathbf{Perf}(k)} && {\mathbf{Perf}(X)}
	\arrow[from=1-1, to=1-3]
	\arrow[from=1-1, to=3-1]
	\arrow["q", dashed, from=1-1, to=3-3]
	\arrow[from=1-3, to=1-5]
	\arrow[from=1-3, to=3-3]
	\arrow["i", from=1-5, to=3-5]
	\arrow["\psi", from=3-1, to=3-3]
	\arrow["s"', from=3-1, to=5-1]
	\arrow["{\ev_1}"', from=3-3, to=3-5]
	\arrow["{\ev_0}"', from=3-3, to=5-3]
	\arrow["{p^\ast}", from=5-1, to=5-3]
\end{tikzcd}
\]
\[
\begin{tikzcd}[sep=small]
	{\widetilde{\mathbf{Syz}}^{\mathcal{L}}(X)} && {\widetilde{\mathbf{Syz}}(X)} && {\mathbf{Perf}(k)^{\Delta^1}} \\
	\\
	{\mathbf{\mathcal{L}}} && {\mathbf{Perf}(X)} && {\mathbf{Perf}(k)}
	\arrow[from=1-1, to=1-3]
	\arrow[from=1-1, to=3-1]
	\arrow["\varphi", from=1-3, to=1-5]
	\arrow[from=1-3, to=3-3]
	\arrow["{\ev_1}", from=1-5, to=3-5]
	\arrow[from=3-1, to=3-3]
	\arrow["{p_\ast}"', from=3-3, to=3-5]
\end{tikzcd}
\]
\end{lemma}

\begin{remark}\label{rem:push__preserves_Perf}
    In the setting of \cref{lemma:equiv_def_Syz_L}, the morphism $p_\ast \colon \mathbf{Perf}(X) \to \mathbf{Perf}(k)$ is well defined by (the proof of) \cref{lem:O-compact}.
\end{remark}

\begin{proof}
	The equivalence with the definition given by the first diagram is obvious.
	For the second diagram, it follows from the functoriality of passing from left to right adjoints, applied to the adjunction $p^\ast \colon \Perf(k) \leftrightarrows
\Perf(X) \colon p_\ast$.
\end{proof}

\begin{definition}\label{def:dSyz_prelim}
    Let $X$ be a proper scheme locally of finite presentation over $k$ and $i \colon \mathcal{L} \to \mathbf{Perf}(X)$ be a morphism in $\dSt_k$.
    We define $\widetilde{\mathbf{Syz}}^{\spl}(X) \in \dSt_k$ via the pullback square in $\dSt_k$
\[
\begin{tikzcd}[sep=small]
	{\widetilde{\mathbf{Syz}}^{\spl}(X)} && {\widetilde{\mathbf{Syz}}(X)} \\
	&& {\mathbf{Perf}(X)^{\Delta^1}} \\
	{\mathbf{Spl}(X)} && {\mathbf{Perf}(X)}
	\arrow["j", from=1-1, to=1-3]
	\arrow[from=1-1, to=3-1]
	\arrow["\psi", from=1-3, to=2-3]
	\arrow["{\ev_1}", from=2-3, to=3-3]
	\arrow[from=3-1, to=3-3]
\end{tikzcd}
\]
    We define similarly $\widetilde{\mathbf{Syz}}^{\mathcal{L}, \spl}(X)$ by replacing $\widetilde{\mathbf{Syz}}(X)$ with $\widetilde{\mathbf{Syz}}^\mathcal{L}(X)$.
\end{definition}

\begin{definition}\label{def:dSyz}
    In the setting of \cref{def:dSyz_prelim}, we define $\mathbf{dSyz}(X) \in \dSt_k$ via the pullback square in $\dSt_k$
\[
\begin{tikzcd}[sep=small]
	{\mathbf{dSyz}(X)} && {\widetilde{\mathbf{Syz}}^{\spl}(X)} \\
	&& {\mathbf{Perf}(X)^{\Delta^1} \times \mathbf{Perf}(k)^{\Delta^1} \times \mathbf{Perf}(k)} \\
	{\mathbf{Perf}^{(0, + \infty)}(X) \times \mathbf{Perf}^{(-\infty, 0)}(k) \times \mathbf{Perf}^{[0,0]}(k)} && {\mathbf{Perf}(X) \times \mathbf{Perf}(k) \times \mathbf{Perf}(k)}
	\arrow[from=1-1, to=1-3]
	\arrow[from=1-1, to=3-1]
	\arrow["(\psi \circ j) \times (\varphi \circ j) \times s", from=1-3, to=2-3]
	\arrow["{\cof \times \fib \times id}", from=2-3, to=3-3]
	\arrow[from=3-1, to=3-3]
\end{tikzcd}
\]
where:
\begin{itemize}
    \item $j \colon \widetilde{\mathbf{Syz}}^{\spl}(X) \to \widetilde{\mathbf{Syz}}(X)$ is the canonical open immersion in \cref{def:dSyz_prelim};
    \item $\varphi \colon \widetilde{\mathbf{Syz}}(X) \to \mathbf{Perf}(k)^{\Delta^1}$, $\psi \colon \widetilde{\mathbf{Syz}}(X) \to \mathbf{Perf}(X)^{\Delta^1}$, and $s \colon \widetilde{\mathbf{Syz}}(X) \to \mathbf{Perf}(k)$ are defined in \cref{lemma:equiv_def_Syz_L};
    \item $\cof \colon \mathbf{Perf}(X)^{\Delta^1} \to \mathbf{Perf}(X)$, $\fib \colon \mathbf{Perf}(k)^{\Delta^1} \to \mathbf{Perf}(k)$ are respectively the cofiber and fiber morphisms;
    \item $\mathbf{Perf}(X)^{(0, + \infty)} \subset \mathbf{Perf}(X)$ and $\mathbf{Perf}(k)^{(-\infty, 0)} \subset \mathbf{Perf}(k)$ are the substacks of perfect complexes with Tor-amplitude contained in $(0, + \infty)$ and in $(-\infty, 0)$.
\end{itemize}
    We define similarly $\mathbf{dSyz}^{\mathcal{L}}(X)$ by replacing $\widetilde{\mathbf{Syz}}^{\spl}(X)$ with $\widetilde{\mathbf{Syz}}^{\mathcal{L}, \spl}(X)$.
    
\end{definition}

\begin{remark}\label{rem:comparison_dSyz}
    Let $\kappa$ be a field over $k$.
    Unwinding the definitions, a $\kappa$-point of $\mathbf{dSyz}(X)$ is given by a morphism $\eta \colon p^\ast W \to F$, where $W \in \Mod_\kappa^\heartsuit$ is a finite dimensional vector space and $F \in \QCoh(X_\kappa)^\heartsuit$ is a perfect simple sheaf, and such that the morphism $\eta$ is surjective and the adjoint morphism $W \to \mathrm{H}^0(p_\ast F)$ is injective.
    \cref{lem:openness_dSyz} shows that this is the set of points of an open substack of $\widetilde{\mathbf{Syz}}^{\spl}(X)$.
    \\ \indent
    Hence, a $\kappa$-point corresponds to a perfect sheaf $F \in \QCoh(X_\kappa)^\heartsuit$ and a vector subspace $W \subset \mathrm{H}^0(X_\kappa,F)$ that globally generates $F$.
    The truncation $\mathbf{Syz}(X) \coloneqq t_0\big(\mathbf{dSyz}(X)\big)$ thus recovers the syzygy stack considered in \cite{FMR25, FMR25b}.
\end{remark}

\begin{lemma}\label{lem:openness_dSyz}
    In the setting of \cref{def:dSyz}, the morphisms
\[
\widetilde{\mathbf{Syz}}^{\spl}(X) \to \widetilde{\mathbf{Syz}}(X), \qquad \mathbf{dSyz}(X) \to \widetilde{\mathbf{Syz}}^{\spl}(X)
\]
and
\[
\widetilde{\mathbf{Syz}}^{\mathcal{L}, \spl}(X) \to \widetilde{\mathbf{Syz}}^{\mathcal{L}}(X), \qquad \mathbf{dSyz}^{\mathcal{L}}(X) \to \widetilde{\mathbf{Syz}}^{\mathcal{L}, \spl}(X)
\]
are representable by open immersions.
    In particular, the derived stacks $\widetilde{\mathbf{Syz}}^{\spl}(X)$ and $\mathbf{dSyz}(X)$ are locally geometric stacks, locally of finite presentation over $k$.
    The same holds for $\widetilde{\mathbf{Syz}}^{\mathcal{L}, \spl}(X)$ and $\mathbf{dSyz}^{\mathcal{L}}(X)$ provided that $\mathcal{L}$ is locally geometric.
\end{lemma}

\begin{proof}
    We treat the case of the first and second morphisms in the statement, the proof for third and fourth morphisms being analogous.
    The first morphism is an open immersion since $\mathbf{Spl}(X) \to \mathbf{Perf}(X)$ is.
    The second morphism is an open immersion since, for any quasi-compact derived scheme $Y$, the morphisms $\mathbf{Perf}(Y)^{(0, + \infty)}, \mathbf{Perf}(Y)^{(-\infty, 0)} \to \mathbf{Perf}(Y)$ are open immersions by \cite[Proposition 6.1.4.5]{Lur18}.
\end{proof}

	As a consequence of \cref{thm:main_DCoh} and \cref{lem:openness_dSyz}, we have:

\begin{corollary}\label{coro:Lagr_dSyz}
	In the setting of \cref{thm:main_DCoh}, the composite morphism in $\dSt_k$
\[
\fib \circ q \colon \mathbf{dSyz}^{\mathcal{L}}(X) \hookrightarrow \widetilde{\mathbf{Syz}}^{\mathbf{\mathcal{L}}}(X) \xrightarrow{\fib \circ q } \mathbf{Perf}(X)
\]
carries a $(1-d)$-shifted Lagrangian structure.
	The same holds for its pulled-back morphism $\fib \circ q \colon \mathbf{dSyz}^{\mathcal{L}}(X) \times_{\mathbf{Perf}(X)} \Spl(X) \to \Spl(X)$ along the open immersion $\Spl(X) \subset \mathbf{Perf}(X)$.
\end{corollary}

\subsection{\texorpdfstring{Shifted symplectic rigidification, after \cite{Park26}}{Shifted symplectic rigidification, after Park--You}}\label{section:shifted_sympl_rig}

	In this paragraph, we recollect the material from \cite{Park26} on shifted symplectic rigidification that we will need later on.
	All the results follow directly from the work of Park--You \cite{Park26} and we include proofs only when the result is not explicitly stated in \textit{loc.cit}.
	The main result that we will need is the shifted Lagrangian rigidification of \cref{thm:rig_Lagr_morph}.\\ \indent
	In this section, following \cite{Park26}, we will assume $k = \mathbb{C}$.
	We fix once and for all a base $B \in \dSt_{\mathbb{C}}$ and an integer $d \leq 0$ and make the standing assumption:
	
\begin{assumption}[{\cite[Assumption 3.2]{Park26}}]
	Assume one of the following two conditions on $B$ and $d$:
\begin{itemize}
	\item $B = \Spec(\mathbb{C})$ and $d\leq 0$;
	\item $B$ is an underived (possibly, higher) Artin stack, locally of finite type over $\mathbb{C}$, and $d \leq -1$.
\end{itemize}
\end{assumption}

\begin{notation}
	To ease the notation, we let 
\[
G \coloneqq B\G_m \ , \qquad \mathfrak{g} \coloneqq \mathbb{T}_{G,id} \ , \qquad T^\ast[d+1]BG \coloneqq T^\ast [d+1](BG/B) \ .
\]
\end{notation}

\begin{recollection}[Symplectic and Hamiltonian $\infty$-categories]\label{recollection:forget_Ham}
	Consider the following $\infty$-categories:
\begin{itemize}
	\item[$Symp_{B,d}^G$:] its objects are derived stacks $\mathcal{X}$ over $BG$ that are $d$-shifted symplectic relative to $BG$. 
	Equivalently, these are derived stacks over $B$ equipped with a $G$-invariant $d$-shifted symplectic structure, the equivalence being given by sending a stack with $G$-action $M$ to the quotient $M / G$.\\
	Morphisms $\mathcal{X}_1 \dashrightarrow \mathcal{X}_2$ are shifted Lagrangian correspondences from $\mathcal{X}_1$ to $\mathcal{X}_2$ over $BG$, that is shifted Lagrangian morphisms $\mathcal{L} \to \mathcal{X}_1 \times_{BG} \mathcal{X}_2$.
	See \cite[Section 2.2]{Park26} for a careful definition.
	\item[$Ham_{B,d}^G$:] its objects are shifted Lagrangian morphisms $\mathcal{L} \to \faktor{\mathfrak{g}^\vee[d]}{G} \simeq T^\ast[d+1] BG$, where $\mathfrak{g} \coloneqq \mathbb{T}_{G, id}$.\\
	Morphisms $\mathcal{L}_1 \dashrightarrow \mathcal{L}_2$ are shifted Lagrangian correspondences over $T^\ast[d+1] BG$, that is, shifted Lagrangian morphisms $\mathcal{L}_3 \to \mathcal{L}_1 \times_{T^\ast[d+1] BG} \mathcal{L}_2$, where the target is $d$-shifted symplectic by the Lagrangian intersection theorem \cite[Theorem 2.9]{PTVV13}.
	See \cite[Section 2.3]{Park26} for a careful definition.
\end{itemize}
	Then, by \cite[Theorem 3.3]{Park26}, the forgetful functor
\begin{equation}\label{eq:equiv_forget_Ham}
\begin{aligned}
\forget & _{Ham} \colon Ham_{B,d}^G \to Symp_{B,d}^G \\
& (\mathcal{L} \xrightarrow{\mu} T^\ast[d+1] BG) \mapsto (\mathcal{L} \xrightarrow{p \circ \mu} \ BG)
\end{aligned}
\end{equation}
is an equivalence, where $p \colon T^\ast[d+1] BG \simeq \faktor{\mathfrak{g}^\vee[d]}{G} \to BG$ is the projection\footnote{The functor $\forget_{Ham}$ is defined as the fiberwise Lagrangian intersection functor along the Lagrangian fibration $p$. 
	Unraveling the definition, this is seen to coincide with the forgetful functor along $p$.
	This description is also explicitly used in the proof of \cite[Lemma 1.9]{Park26}.}.
	Notice that $\mathcal{L}$ is $d$-shifted symplectic over $BG$ by \cite[Lemma 1.10]{Saf23} and \cite[Notation 2.5]{Park26}.
	Notice also that $\forget_{Ham}$ sends a Lagrangian morphism over $T^\ast[d+1] BG$
\[
\mathcal{L}_3 \to \mathcal{L}_1 \times_{T^\ast[d+1] BG} \mathcal{L}_2 \ ,
\]
that is a morphism in $Ham_{B,d}^G$, to the Lagrangian morphism over $BG$ (\cite[Lemma 2.4.4]{Par24})
\[
\mathcal{L}_3 \to \mathcal{L}_1 \times_{BG} \mathcal{L}_2 \ ,
\]
that is a morphism in $Symp_{B,d}^G$.
\end{recollection}

\begin{definition}\label{def:symplectic_rigidification_Ham}
    In the setting of \cref{recollection:forget_Ham}, let $\mu \colon \mathcal{L} \to T^\ast[d+1] BG$ be a shifted Lagrangian morphism, i.e., an object of $Ham_{B,d}^G$.
    The shifted symplectic rigidification $\mathcal{L} /\!\!/ B \G_m$ of $\mu$ is the $d$-shifted symplectic stack defined as the Lagrangian intersection (\cite[Theorem 2.9]{PTVV13})
    \begin{align}\label{diag:shifted_symplectic_rigidification}
        \begin{tikzcd}[ampersand replacement = \&]
            \mathcal{L} /\!\!/ B \G_m \ar[r, "\lambda"] \ar[d, "\psi"'] \& BG \ar[d, "0"] \\
            \mathcal{L} \ar[r, "\mu"'] \& T^\ast[d+1]BG \ .
        \end{tikzcd}
    \end{align}
\end{definition}

\begin{recollection}\label{recollection:symplectic_action}
	In the setting of \cref{recollection:forget_Ham}, let $M \in \dSt_B$ be a locally geometric stack carrying a $G$-action and a $d$-shifted symplectic structure.
	The $d$-shifted symplectic form on $M$ has a unique (up to contractible space of choices) $G$-invariant lift by \cite[Corollary 3.10]{Park26}.
	Equivalently, by \cref{recollection:forget_Ham}, we have that $M / G \in Symp_{B,d}^{G}$.
	Through the equivalence \eqref{eq:equiv_forget_Ham}, there is an associated $G$-equivariant moment map $\mu_M \colon M \to \mathfrak{g}^\vee[d]$ and an associated Lagrangian structure on the quotient morphism
\[
\mu_M / G \colon M/G \to T^\ast[d+1]BG \simeq \faktor{\mathfrak{g}^\vee[d]}{G} \ .
\]
	Moreover, by \cite[Proposition 3.1]{Park26}, given any shifted Lagrangian correspondence
\[
\begin{tikzcd}[sep=small]
	&& {\mathcal{L}} && \\
	{M_1} &&&& {M_2}
	\arrow["{i_1}"', from=1-3, to=2-1]
	\arrow["{i_2}", from=1-3, to=2-5]
\end{tikzcd}
\]
where $\mathcal{L}, M_1, M_2$ carry $G$-actions and $i_1, i_2$ are $G$-equivariant, there exists a unique (up to contractible spaces of choice) lift of the shifted Lagrangian correspondence to a $G$-equivariant shifted Lagrangian correspondence.
	Hence, we get a Lagrangian correspondence over $BG$
\begin{equation}\label{eq:Lagr_corr_forget}
\begin{tikzcd}[sep=small]
	&& {\mathcal{L} /G} && \\
	{M_1 /G} &&&& {M_2 /G}
	\arrow["{i_1 /G}"', from=1-3, to=2-1]
	\arrow["{i_2 /G}", from=1-3, to=2-5]
\end{tikzcd}
\end{equation}
	Through the equivalence \eqref{eq:equiv_forget_Ham}, there exist unique (up to contractible spaces of choices) shifted Lagrangian morphisms
\begin{align*}
\mu_1 \colon M_1 / G \to T^\ast[d+1] BG \ , \\
\mu_2 \colon M_2 / G \to T^\ast[d+1] BG \ , \\
\overline{i} \colon \mathcal{L} / G \to  \big( M_1 / G \big)  \times_{T^\ast[d+1] BG}  \big( M_2 / G \big)
\end{align*}
such that $\forget_{Ham}(\mu_i)$ is the structural morphism $M_i / G \to BG$ and $\forget_{Ham}(\overline{i})$ is the Lagrangian correspondence \eqref{eq:Lagr_corr_forget}.
\end{recollection}

\begin{definition}\label{def:symplectic_rigidification}
	In the setting of \cref{recollection:symplectic_action}, we define the shifted symplectic rigidification $M /\!\!/ G$ of $M$ as the shifted symplectic rigidification of $\mu_M / G$ in the sense of \cref{def:symplectic_rigidification_Ham}.
\end{definition}

\begin{remark}\label{rem:truncation_of_rigidification}
    As taking the classical truncation of higher Artin stacks commutes with limits, we see from \eqref{diag:shifted_symplectic_rigidification} that the symplectic rigidification $M /\!\!/ G$ and the ordinary rigidification $M / G$ have the same classical truncation, when considered as higher Artin stacks.
    It then follows that $t_0(M /\!\!/ G) \simeq t_0(M / G)$.
\end{remark}




\begin{theorem}[Shifted Lagrangian rigidification]\label{thm:rig_Lagr_morph}
	In the setting of \cref{def:symplectic_rigidification}, let $i \colon \mathcal{L} \to M$ be a $G$-equivariant morphism carrying a shifted Lagrangian structure.
	Then, the induced morphism
\[
\overline{i} \colon \mathcal{L} / G \to M /\!\!/ G
\]	
from the rigidification of $\mathcal{L}$ to the shifted symplectic rigidification of $M$ carries a shifted Lagrangian structure. 
\end{theorem}

\begin{proof}
	The shifted Lagrangian morphism $i$ defines a Lagrangian correspondence $B \dashrightarrow M$ where the morphisms are $G$-equivariant.
	By \cref{recollection:symplectic_action}, there exist unique (up to contractible spaces of choices) shifted Lagrangian morphisms
\begin{align*}
\mu_1 \colon BG \to T^\ast[d+1] BG \ , \\
\mu_2 \colon M / G \to T^\ast[d+1] BG \ , \\
\overline{i} \colon \mathcal{L} / G \to BG \times_{T^\ast[d+1] BG}  (M / G )
\end{align*}
such that $\forget_{Ham}(\mu_1) \simeq id_{BG}$, $\forget_{Ham}(\mu_2) = M / G \to BG$ is equivalent to the canonical projection, and $\forget_{Ham}(\overline{i})$ is the shifted Lagrangian correspondence $BG \dashrightarrow M/G$ associated to the shifted Lagrangian morphism $i / G$ over $BG$.
	Hence, we are left to identify $\mu_1 \colon BG \to T^\ast[d+1] BG$ with the zero section.
	The latter claim follows directly from the fact that $\forget_{Ham}(BG \xrightarrow{0} T^\ast[d+1] BG) = (BG \xrightarrow{id} BG)$, where on the right $BG$ is considered with its trivial $d$-shifted symplectic structure relative to $BG$ itself.
\end{proof}

\begin{remark}\label{rem:PY_shifted_Lagr_rig}
    As communicated to us by Park and You, \cref{thm:rig_Lagr_morph} also follows by the adjunction in \cite[Corollary 3.5]{Park26}.
    More precisely, under the equivalence $\Map(M, \mathrm{triv}(B)) \simeq \Map(M /\!\!/ G, B)$, the shifted Lagrangian correspondence $M \dashrightarrow \mathrm{triv}(B)$ given by $i$ corresponds to the shifted Lagrangian correspondence $M /\!\!/ G \dashrightarrow B$ given by the $\overline{i}$ appearing in \cref{thm:rig_Lagr_morph}.
\end{remark}

	We now apply the precedent discussion to the moduli of perfect complexes on a Calabi--Yau scheme.
	Notice that we will now denote $d \coloneqq \dim X$, so that the shift of the symplectic structure on $\mathbf{Perf}(X)$ is $2-d$.

\begin{recollection}\label{BG_m_action_Perf}
	Let $X$ be a scheme over $k$.
	Then, by \cite[Construction 4.7]{Park26}, there is a natural $B\G_m$-action on $\mathbf{Perf}(X)$ given by the constant automorphisms of perfect complexes.
	One immediately checks that the above $B\G_m$-action restricts to a $B\G_m$-action on $\mathbf{Perf}^{[0,0]}(X)$ and $\Spl(X)$.
\end{recollection}

\begin{definition}
	Let $X$ be a proper Calabi--Yau scheme over $\mathbb{C}$ of dimension $d \geq 2$.
    We let $\mathbf{Perf}(X) /\!\!/ B \G_m$ be the symplectic rigidification of $\mathbf{Perf}(X)$ with respect to the $B\G_m$-action of \cref{BG_m_action_Perf}, in the sense of \cref{def:symplectic_rigidification}.
    We define $\mathbf{Perf}^{[0,0]}(X) /\!\!/ B\G_m $ and $\Spl(X) /\!\!/ B\G_m$ similarly.
\end{definition}

\begin{corollary}\label{coro:Lagr_rig_pullback}
	Let $X$ be a strong Calabi--Yau scheme over $\mathbb{C}$ of dimension $d \geq 2$.
	Then, the morphism $p^\ast \colon \mathbf{Perf}(\mathbb{C}) \to \mathbf{Perf}(X)$ induces a shifted Lagrangian morphism
\[
p^\ast \colon \mathbf{Perf}(\mathbb{C}) / B \G_m \to \mathbf{Perf}(X) /\!\!/ B \G_m \ .
\]
	The same holds by replacing $\mathbf{Perf}(-)$ with $\mathbf{Perf}^{sim}(-)$ and $\Spl(-)$.
\end{corollary}

\begin{proof}
	It follows directly from \cref{coro:Lagr_strong} and \cref{thm:rig_Lagr_morph}.
\end{proof}

\begin{corollary}\label{coro:rigid_Lagr_Spl}
	In the setting of \cref{thm:main_DCoh}, let $k = \mathbb{C}$ and $d \geq 2$.
	If $\mathbf{\mathcal{L}}$ carries a $B\G_m$-action and $i$ is $B\G_m$-equivariant, then the rigidified morphism
\[
(\fib \circ q)^{rig} \colon \widetilde{\mathbf{Syz}}^{\mathbf{\mathcal{L}}}(X) / B\G_m \to \mathbf{Perf}(X) /\!\!/ B\G_m
\]
carries a $(1-d)$-shifted Lagrangian structure.
	The same holds for
\[
(\fib \circ q)^{rig} \colon \big(\mathbf{dSyz}^{\mathcal{L}}(X) \times_{\mathbf{Perf}(X)} \Spl(X)\big) / B\G_m\to \Spl(X)  /\!\!/ B\G_m \ .
\]
\end{corollary}

\begin{proof}
	Since all the morphisms in \eqref{eq:def_Syz_tilde} are $B\G_m$-equivariant, it follows that $\widetilde{\mathbf{Syz}}^{\mathbf{\mathcal{L}}}(X)$ carries a $B\G_m$-action and that $\fib \circ q$ is $B\G_m$-equivariant.
	The first claim then follows by \cref{thm:main_DCoh} and \cref{thm:rig_Lagr_morph}.
	The second claim is proven analogously replacing \cref{thm:main_DCoh} with \cref{coro:Lagr_dSyz}.
\end{proof}
	
	For future reference, we also record the following result of Park--You \cite{Park26} on the rigidification of the Lagrangian correspondence of \cref{prop:Lagrangian_ev_fib}.
	
\begin{recollection}\label{recollection:rigid_correspondece}
	Let $X$ be a proper Calabi--Yau scheme over $\mathbb{C}$ of dimension $d \geq 2$.
	Then, by \cref{prop:Lagrangian_ev_fib} we have a Lagrangian correspondence
\[
\begin{tikzcd}[sep=small]
	& {\mathbf{Perf}(X)^{\Delta^1}} && \\
	{\mathbf{Perf}(X) \times_\mathbb{C} \mathbf{Perf}(X) } &&& {\mathbf{Perf}(X)}
	\arrow["{(\fib, \ev_1)}"', from=1-2, to=2-1]
	\arrow["{\ev_0}", from=1-2, to=2-4]
\end{tikzcd}
\]
	By \cite[Theorem D]{Park26} (see also \cite[Construction 4.18]{Park26}), there exists a derived stack $\mathbf{Perf}(X)^{\Delta^1, rig}$ such that $t_0 \mathbf{Perf}(X)^{\Delta^1, rig} \simeq t_0  \big(\mathbf{Perf}(X)^{\Delta^1} /\!\!/ B\G_m\big)$ and such that there is an induced Lagrangian correspondence
\[
\begin{tikzcd}[sep=small]
	& {\mathbf{Perf}(X)^{\Delta^1, rig}} && \\
	{\big(\mathbf{Perf}(X) /\!\!/ B\G_m\big) \times_\mathbb{C} \big(\mathbf{Perf}(X) /\!\!/ B\G_m\big) } &&& {\mathbf{Perf}(X) /\!\!/ B\G_m}
	\arrow["{(\fib, \ev_1)}"', from=1-2, to=2-1]
	\arrow["{\ev_0}", from=1-2, to=2-4]
\end{tikzcd}
\]
\end{recollection}

\subsection{From derived algebraic stacks to classical algebraic spaces}

	In this paragraph, we use \cref{thm:rig_Lagr_morph} to descend the shifted Lagrangian morphism of \cref{coro:Lagr_strong} to a Lagrangian morphism of classical algebraic spaces of simple perfect complexes.
    This is the key input to recover and generalize the results of Fantechi--Miró-Roig \cite{FMR25, FMR25b}.\\ \indent
    We will need the following result of Bozec--Calaque--Scherotzke \cite{BCS24b}, whose proof strategy they attribute to Benjamin Hennion.
    The result in \textit{loc.cit.} is stated for smooth algebraic varieties, but the proof applies to smooth algebraic spaces with no change.

\begin{proposition}[{\cite[Proposition 4.26]{BCS24b}}]\label{prop:classical_lagr_from_derived}
    Consider a Lagrangian morphism $g \colon \mathcal{L} \to Y$, with $Y$ a classical smooth and symplectic algebraic space over $k$.
    Denote by $t_0(\mathcal{L}) \to \mathcal{L}$ the truncation of $\mathcal{L}$, and let $Z \subset Y$ be the schematic image of $t_0(\mathcal{L}) \to Y$.
    Then $Z \subset Y$ is a classical Lagrangian sub-algebraic space.
\end{proposition}

\begin{recollection}\label{recollection:gerbes_simple}
    Let $X$ be a smooth and proper Calabi--Yau surface.
    The substacks $\mathbf{Perf}^{sim}(X), \mathbf{Perf}(X)^{[0,0]}, \Spl(X) \subset \mathbf{Perf}(X)$ being open substacks of a $0$-shifted symplectic stack (\cref{coro:sympl_Calbi_Yau}), are themselves $0$-shifted symplectic.
    As explained in \cite[Section 3.1]{PTVV13}, the truncation $t_0(\mathbf{Perf}(X)^{sim})$ is a $\G_m$-gerbe over a smooth algebraic space $\underline{Perf}^{sim}(X)$ (smoothness is proven in \cite[Theorem 3.3]{Ina11}).
    Moreover, the algebraic space $\underline{Perf}^{sim}(X)$ inherits a symplectic form from $\mathbf{Perf}(X)^{sim}$, whose underlying pairing on the tangent space at a simple perfect complex $F$ is given by
\[
\mathrm{Ext}^1(F,F) \times \mathrm{Ext}^1(F,F) \xrightarrow{\cap} \mathrm{Ext}^2(F,F) \xrightarrow{tr} \Ext^2(\OO_X, \OO_X) \simeq k \ .
\]
    We will provide a generalization of the above discussion in \cref{prop:smooth_rigidification}.
\end{recollection}

    As already noticed in \cite[Section 3.1]{PTVV13}, in the setting of \cref{recollection:gerbes_simple}, the construction of a non-degenerate $2$-form $\omega$ on $\underline{Perf}^{sim}(X)$ via classical tools is easy, but the proof that $\omega$ is closed requires heavy computations (see \cite{Ina11}).
    Similarly, in \cite{FMR25b} the authors prove the existence of a non-degenerate $2$-form $\omega$, but the closedness of $\omega$ is only proven under a smoothability assumption.
    The approach of \cite{PTVV13} described in \cref{recollection:gerbes_simple} uses derived algebraic geometry to construct a \textit{closed} $2$-form on $\underline{Perf}^{sim}(X)$, which is readily verified to be non-degenerate.
    We will take a similar path to \cite{PTVV13}, but employing different tools: we will consider the \textit{shifted symplectic rigidification} of $\Spl(X)$, which we then prove to be equivalent to the classical rigidification $\underline{Spl}(X)$ of $t_0(\Spl(X))$ when $\dim X = 2$.

\begin{definition}\label{def:classical_rigidifiaction_Spl}
    Let $X$ be a proper scheme locally of finite presentation over $k$.
    We define $\underline{Spl}(X) \coloneqq t_0 \big(\mathbf{Spl}(X)\big) / B\G_m$ as the algebraic space over $k$ obtained as the classical rigidification of the classical stack $t_0 \mathbf{Spl}(X)$.
\end{definition}

\begin{remark}\label{rem:classica_vs_derived_rig}
	In the setting of \cref{def:classical_rigidifiaction_Spl}, notice that we have a canonical equivalence $\underline{Spl}(X) \coloneqq t_0 \big(\mathbf{Spl}(X)\big) / B\G_m \simeq t_0 \big(\mathbf{Spl}(X) / B\G_m \big)$.
\end{remark}

    The following result is certainly known to Park--You, at least in the case where $X$ is smooth, see \cite[Theorem B]{Park26}.
    Since we could not locate a proof of it, we include it here for completeness.

\begin{proposition}\label{prop:smooth_rigidification}
    Let $X$ be a proper Calabi--Yau scheme over $\mathbb{C}$ of dimension $2$.
    The derived symplectic rigidification $\mathbf{Perf}^{sim}(X) /\!\!/ B \G_m$ and $\Spl(X) /\!\!/ B\G_m$ are classical smooth algebraic spaces over $\mathbb{C}$.
    In particular, we have an equivalence $\mathbf{Spl}(X) /\!\!/ B \G_m \simeq \underline{Spl}(X)$.
\end{proposition}

    For the next proof, we recall that a morphism $f \colon X \to Y$ of derived stacks is \emph{locally almost of finite presentation} if the following condition holds: given $n \geq 0$, for every filtered diagram $(A_i)_{i \in I}$ of $n$-truncated algebras with colimit $A$, then the canonical morphism
\[
\colim_{i \in I}X(A_i) \to X(A) \times_{Y(A)} \big(\colim_{i \in I} Y(A_i)\big)
\]
is an equivalence, see \cite[Definition 17.4.1.1-(b)]{Lur18}.

\begin{proof}
    By definition and \cref{rem:truncation_of_rigidification}, both derived Artin stacks $\mathbf{Perf}^{sim}(X) /\!\!/ B \G_m$ and $\Spl(X) /\!\!/ B\G_m$ are $0$-shifted symplectic algebraic spaces.
    Moreover, $\Spl(X) /\!\!/ B\G_m$ is an open sub-algebraic space of $\mathbf{Perf}^{sim}(X) /\!\!/ B \G_m$. Indeed, this can be checked locally in the \'etale topology, so it follows from the fact that $\Spl(X) \to \mathbf{Perf}^{sim}(X)$ is an open immersion.
    Hence, it is enough to show that the cotangent complex of $\mathbf{Perf}^{sim}(X) /\!\!/ B \G_m$ has Tor-amplitude contained in $[0,0]$.
    Since $\mathbf{Perf}^{sim}(X) /\!\!/ B \G_m$ has perfect cotangent complex, being a pullback of locally almost of finite presentation derived stacks over $k$, we can equivalently show that the tangent complex of $\mathbf{Perf}^{sim}(X) /\!\!/ B \G_m$ has Tor-amplitude contained in $[0,0]$.
    Let $\mathbb{C} \subset \kappa$ be a field.
    By \cite[Corollary 6.1.4.7]{Lur18}, we are left to show that for every $\kappa$-point $x_F \colon \Spec(\kappa) \to \mathbf{Perf}^{sim}(X)$, corresponding to a simple perfect complex $F$ on $X_\kappa$, we have that $\mathbb{T}_{\mathbf{Perf}^{sim}(X)/\!\!/ B \G_m, F}$ is concentrated in degree $0$.
    More precisely, we will show that
\[
\mathbb{T}_{\mathbf{Perf}^{sim}(X) /\!\!/ B \G_m, F} \simeq \mathrm{H}^0\mathbb{T}_{\mathbf{Perf}^{sim}(X), F} \simeq \mathrm{Ext}^1(F,F) \footnote{as the $\kappa$-points of $\mathbf{Perf}^{sim}(X), \mathbf{Perf}^{sim}(X) / B\G_m$ and $\mathbf{Perf}^{sim}(X) /\!\!/ B\G_m$ are the same, we identify them in the formulas involving tangent complexes.
    In particular, if $\chi$ denotes the morphism $\chi \colon \mathbf{Perf}^{sim}(X) \to \mathbf{Perf}^{sim}(X) /\!\!/ B\G_m$, by $\mathbb{T}_{\mathbf{Perf}^{sim}(X) /\!\!/ B\G_m, F}$ we mean $x_F^\ast \chi^\ast \mathbb{T}_{\mathbf{Perf}^{sim}(X) /\!\!/ B\G_m}$, and similarly for $\lambda \colon \mathbf{Perf}^{sim}(X) /\!\!/ B \G_m \to \mathbf{Perf}^{sim}(X) / B\G_m$.}\ .
\]
    Notice that
\[
\mathbb{T}_{\mathbf{Perf}^{sim}(X) / B \G_m, F} \simeq \cof\left(\kappa[1] \xrightarrow{g} \mathbb{T}_{\mathbf{Perf}^{sim}(X), F}\right) \simeq \tau_{\leq 0} \mathbb{T}_{\mathbf{Perf}^{sim}(X), F}  \ ,
\]
as $g$ is identified with the inclusion of the scalar multiplications morphisms of the simple complex $F$, under the equivalence $\mathbb{T}_{\mathbf{Perf}^{sim}(X), F} \simeq \Hom_{\QCoh(X_\kappa)}(F,F)[1]$.
    Moreover,  since the zero section $0: B^2 \G_m \to T^\ast [1]B^2 \G_m$ is $0$-shifted Lagrangian by \cite[Corollary 2.4]{Cal19}, we have an equivalence
\[
\mathbb{T}_0 \simeq \mathbb{L}_{B^2 \G_m} \simeq \OO_{B^2\G_m} \otimes_k \mathfrak{g}_a[-2] \ ,
\]
where $\mathfrak{g}_a$ is the Lie algebra of $\G_m$.
    The pullback square in \cref{def:symplectic_rigidification} then gives a diagram in $\Mod_k$ with pushout and pullback squares
\[
\begin{tikzcd}[sep=small]
	{\mathbb{T}_{\mathbf{Perf}^{sim}(X) /\!\!/ B \G_m, F}} && {x_F^\ast \lambda^\ast\mathbb{T}_{B^2\G_m}} && 0 \\
	\\
	{\tau_{\leq 0} \mathbb{T}_{\mathbf{Perf}^{sim}(X), F}} && {x_F^\ast \lambda^\ast 0^\ast\mathbb{T}_{\mathbb{T}^\ast[1] B^2\G_m}} && {(\mathfrak{g}_a \otimes_k \kappa[-2])[1] \simeq \mathfrak{g}_a \otimes_k \kappa[-1]}
	\arrow[from=1-1, to=1-3]
	\arrow[from=1-1, to=3-1]
	\arrow[from=1-3, to=1-5]
	\arrow[from=1-3, to=3-3]
	\arrow[from=1-5, to=3-5]
	\arrow[from=3-1, to=3-3]
	\arrow[from=3-3, to=3-5]
\end{tikzcd}
\]
	The exterior pullback and pushout square shows that $\mathbb{T}_{\mathbf{Perf}^{sim}(X) /\!\!/ B \G_m, F}$ is coconnective.
	Since $\mathbf{Perf}^{sim}(X) /\!\!/ B \G_m$ is $0$-shifted symplectic by construction (\cref{def:symplectic_rigidification}), it follows by autoduality that $\mathbb{T}_{\mathbf{Perf}^{sim}(X) /\!\!/ B \G_m, F}$ is also connective, hence proving the first claim of \cref{prop:smooth_rigidification}.\\ \indent
	Finally, by \cref{rem:truncation_of_rigidification} and \cref{rem:classica_vs_derived_rig} and since $\mathbf{Spl}(X) /\!\!/ B \G_m$ is underived, we have equivalences
\[
\mathbf{Spl}(X) /\!\!/ B \G_m \simeq t_0\big(\mathbf{Spl}(X) /\!\!/ B \G_m\big) \simeq t_0\big(\mathbf{Spl}(X) / B \G_m\big) \simeq t_0\big(\mathbf{Spl}(X)\big)  / B \G_m \eqqcolon \underline{Spl}(X) \ .
\]
\end{proof}

\begin{corollary}\label{coro:classic_Lagr_from_divisors}
	In the setting of \cref{thm:CY_divisors}, let $k = \mathbb{C}$, $d=2$ and assume that $X$ is strong Calabi--Yau.
	Consider the morphism
\[
i^!_{spl} \colon \mathbf{Perf}_{prop}(Y) \times_{\mathbf{Perf}(X)} \Spl(X) \to \Spl(X)
\]
obtained by pulling back $i^! \colon \mathbf{Perf}_{prop}(Y) \to \mathbf{Perf}(X)$ along $\Spl(X) \to \mathbf{Perf}(X)$.
	Then $i^!_{spl}$ carries a shifted Lagrangian structure and the schematic image $L$ of its rigidification
\[
t_0\big(\mathbf{Perf}_{prop}(Y) \times_{\mathbf{Perf}(X)} \Spl(X)\big) / B\G_m \to \underline{Spl}(X)
\]
	is a classical Lagrangian sub-algebraic space of $\underline{Spl}(X)$.
\end{corollary}

\begin{proof}
	The morphism $i^!_{spl}$ has a shifted Lagrangian structure induced by \cref{thm:Lagr_divisors}, since $\Spl(X) \to \mathbf{Perf}(X)$ is an open immersion.
	The last claim then follows by combining \cref{thm:rig_Lagr_morph} with \cref{prop:smooth_rigidification}.
\end{proof}
    
    If $X$ is a K3 surface, then the structure sheaf $\OO_X$ is a Lagrangian sub-algebraic space of $\underline{Spl}(X)$.
    This can be proven directly, as $\OO_X$ is an isolated point of $\underline{Spl}(X)$ since $\mathrm{Ext}^1(\OO_X, \OO_X) \simeq H^1(X, \OO_X) = 0$ (see e.g., the proof of \cite{FMR25} for $L = \OO_X$).
    One can prove the same result for (possibly singular) Calabi--Yau schemes, by replacing $\underline{Spl}(X)$ with $\mathbf{Spl}(X) /\!\!/ B\G_m$ if $d > 2$.

\begin{corollary}\label{coro:classical_Lagrangian_Spl}
    Let $X$ be a strong Calabi--Yau scheme over $\mathbb{C}$ of dimension $2$ and $p \colon X \to \Spec(\mathbb{C})$ be the structure morphism.
    Let $L$ be the schematic image of
\[
p^\ast \colon \Spec(\mathbb{C}) \simeq \mathbf{Spl}(\mathbb{C}) /\ B\G_m \to \mathbf{Spl}(X) /\!\!/ B\G_m \ .
\]
	Then $L$ is a classical Lagrangian sub-algebraic space of $\mathbf{Spl}(X) /\!\!/ B\G_m \simeq \underline{Spl}(X)$ and $L$ is canonically identified with the point $\OO_X$.
\end{corollary}

\begin{proof}
    By \cref{coro:Lagr_rig_pullback}, \cref{prop:classical_lagr_from_derived} and \cref{prop:smooth_rigidification}, we have that $L$ is a classical algebraic space of $\underline{Spl}(X)$.
    Notice that $\mathbf{Spl}(\mathbb{C}) \simeq B\G_m$, so that $\Spec(\mathbb{C}) \simeq \mathbf{Spl}(\mathbb{C}) /\ B\G_m$.
    To prove the last claim, notice also that  $\OO_X$ is an isolated point of $\underline{Spl}(X)$ since $\mathrm{Ext}^1(\OO_X, \OO_X) \simeq H^1(X, \OO_X) = 0$.
    Hence, the rigidified morphism  $p^\ast \colon \Spec(\mathbb{C}) \to \underline{Spl}(X)$ is the (open and) closed immersion corresponding to $\OO_X$.
\end{proof}

	By combining \cref{coro:rigid_Lagr_Spl}, \cref{prop:classical_lagr_from_derived}, and \cref{prop:smooth_rigidification}, we get:

\begin{corollary}\label{coro:classic_Lagr_from_derived_Lagr}
	In the setting of \cref{coro:rigid_Lagr_Spl} let $d=2$. Then, the schematic image of the morphism
\[
t_0(\fib \circ q)^{rig} \colon t_0\big(\mathbf{dSyz}^{\mathcal{L}}(X) \times_{\mathbf{Perf}(X)} \Spl(X)\big) / B\G_m\to \underline{Spl}(X)
\]
is a classical Lagrangian sub-algebraic space.
\end{corollary}

	\cref{coro:classic_Lagr_from_derived_Lagr} is a first generalization of \cite[Theorem 1.3]{FMR25b}, where the input is a shifted Lagrangian morphism $\mathcal{L} \to \mathbf{Perf}(X)$.
	In the last part of this paragraph we briefly explain how to obtain the same result by starting with a classical Lagrangian sub-algebraic space $L \to \underline{Spl}(X)$.

\begin{definition}\label{recollection:classical_results}
	Let $X$ be a proper Calabi--Yau scheme over $\mathbb{C}$ of dimension $d \geq 2$ and $i \colon L \to \Spl(X) /\!\!/ B\G_m $ be a shifted Lagrangian morphism.
	We define $\widetilde{\mathbf{Syz}}(X)^{L,rig} \in \dSt_{\mathbb{C}}$ via the diagram with pullback square in $\dSt_{\mathbb{C}}$
\[
\begin{tikzcd}[sep=small]
	{\widetilde{\mathbf{Syz}}(X)^{L,rig}} && {\big(\mathbf{Perf}(\mathbb{C}) / B\G_m\big) \times L} \\
	\\
	{\mathbf{Perf}(X)^{\Delta^1, rig}} && {\big(\mathbf{Perf}(X) /\!\!/ B\G_m\big) \times_\mathbb{C} \big(\mathbf{Perf}(X) /\!\!/ B\G_m\big)} \\
	\\
	{\mathbf{Perf}(X) /\!\!/ B\G_m}
	\arrow["{(\widetilde{\ev}_0, \widetilde{\ev}_1)}", from=1-1, to=1-3]
	\arrow["q"', from=1-1, to=3-1]
	\arrow["{(p^\ast, j_X \circ i)}", from=1-3, to=3-3]
	\arrow["{(\ev_0, \ev_1)}"', from=3-1, to=3-3]
	\arrow["\fib"', from=3-1, to=5-1]
\end{tikzcd}
\]
where $j_X \colon \Spl(X) /\!\!/ B\G_m \to \mathbf{Perf}(X) /\!\!/ B\G_m$ is the open immersion induced by the open immersion $\Spl(X) \to \mathbf{Perf}(X)$.
\end{definition}

\begin{definition}\label{def:Spl_rig}
	In the setting of \cref{recollection:classical_results}, we define $\mathbf{Syz}(X)^{L,rig} \in \dSt_{\mathbb{C}}$ via the pullback square in $\dSt_{\mathbb{C}}$
\[
\begin{tikzcd}[sep=small]
	{\mathbf{Syz}(X)^{L,rig}} && {\widetilde{\mathbf{Syz}}(X)^{L,rig}} \\
	\\
	{\Spl(X) /\!\!/ B\G_m} && {\mathbf{Perf}(X) /\!\!/ B\G_m}
	\arrow["j", from=1-1, to=1-3]
	\arrow["{\fib \circ q}"', from=1-1, to=3-1]
	\arrow["{\fib \circ q}", from=1-3, to=3-3]
	\arrow["{j_X}"', from=3-1, to=3-3]
\end{tikzcd}
\]
\end{definition}

	By \cref{recollection:classical_results}, the following is a generalization of \cite[Theorem 1.3-(2)]{FMR25b}.
	While the statement is purely in terms of classical algebraic geometry, the proof needs all the flexibility of derived algebraic geometry, as the proofs presented above show.
	
\begin{theorem}\label{thm:derived_FMR}
	In the setting of \cref{recollection:classical_results}, assume furthermore that $X$ is strong Calabi--Yau.
	Then, the morphisms
\[
\fib \circ q \colon \widetilde{\mathbf{Syz}}(X)^{L,rig} \to \mathbf{Perf}(X) /\!\!/ B\G_m \ , \qquad \fib \circ q \colon \mathbf{Syz}(X)^{L,rig} \to \Spl(X) /\!\!/ B\G_m
\]
are $(1-d)$-shifted Lagrangian.
\end{theorem}

\begin{proof}
	The Lagrangianity of the second morphism in the claim follows directly from the Lagrangianity of the first morphism and the fact that $j$ and $j_X$ in \cref{def:Spl_rig} are open immersions. 
	For the first morphism, the result follows by \cref{thm:composition_lagr_morph}, \cref{coro:Lagr_rig_pullback} (together with the fact that $\mathbf{Perf}^{[0,0]}(X) \subset \mathbf{Perf}(X)$ is an open substack), \cref{thm:rig_Lagr_morph}, and \cref{recollection:rigid_correspondece}.
\end{proof}

\begin{corollary}\label{thm:classical_results}
	In the setting of \cref{thm:derived_FMR}, assume furthermore that $d=2$.
	Then, the schematic image of
\begin{equation}\label{eq:schemtaic_image}
t_0(\fib \circ q) \colon t_0 \mathbf{Syz}(X)^{L,rig} \to \underline{Spl}(X)
\end{equation}
is a classical Lagrangian sub-algebraic space.
\end{corollary}

\begin{proof}
	It follows directly by \cref{prop:classical_lagr_from_derived}, \cref{prop:smooth_rigidification}, and \cref{thm:derived_FMR}.
\end{proof}

\begin{example}
	In the setting of \cref{thm:classical_results}, we can take $L$ to be the point of $\underline{Spl}(X)$ representing $\OO_X$ by \cref{coro:classical_Lagrangian_Spl} (or any other isolated point).
\end{example}

\begin{example}
	In the setting of \cref{thm:classical_results}, we can take $L$ to be as in \cref{coro:classic_Lagr_from_divisors} if we further assume that $X$ is a K3 surface.
	See \cite{Bea04} for a characterization of K3 surfaces arising as anti-canonical divisor of a Fano threefold.
\end{example}

\begin{remark}\label{rem:FMR_comparison}
	\cref{thm:classical_results} is stable under passing to open substacks of the source $\mathbf{Syz}(X)^{L,rig}$.
	Hence, the classical algebraic space $\mathcal{G}^0_L$ considered in \cite{FMR25, FMR25b} is canonically identified with the classical open sub-algebraic space of $t_0 \mathbf{Syz}(X)^{L,rig}$, obtained by imposing the same conditions appearing in \cref{rem:comparison_dSyz} and the vanishing of the $\mathrm{H}^{1}$ on the families of perfect sheaves that we are parametrizing.
	The induced morphism $\fib \circ q \colon \mathcal{G}^0_L \to \underline{Spl}(X)$ is a locally closed embedding as shown in \cite{FMR24}.
\end{remark}

\begin{remark}
    By \cref{rem:FMR_comparison}, if $X$ is a smoothable $K3$-surface, \cref{thm:classical_results} recovers the main result of \cite{FMR25b} up to passing to an open substack of $\mathcal{G}^0_L \subset \mathbf{Syz}(X)^{L,rig}$ of \cref{def:Spl_rig}.
    Moreover, we notice that if we fix a connected component $\underline{Spl}(X)_\alpha$ of $\underline{Spl}(X)$ and denote by $L_\alpha \coloneqq \underline{Spl}(X)_\alpha \times_{\underline{Spl}(X)} L$, the newly constructed Lagrangian is a Lagrangian sub-algebraic space of $\underline{Spl}(X)_\beta$, where $\beta \in \pi_0(\underline{Spl}(X))$ is a different connected component.
    The precise relation between $\alpha$ and $\beta$, formulated in terms of Chern classes, can be found in \cite[Theorem 1.1]{FMR25}.
\end{remark}

\appendix

\section{\texorpdfstring{The proof of \cref{thm:comparison_Perf_corpus}}{The proof of thm:comparison-Perf-corpus}}\label{Appendix_FM}

\subsection{Fourier--Mukai transforms}
	We start by giving a description of the $A$-points of $\M_{\IndCoh(X)}$ for $A$ a regular ring.
    This will be used to prove \cref{thm:comparison_moduli_Perf_DCoh}.

\begin{definition}[{\cite[Definition 1.2.3]{BZNP17}}]\label{notation:Coh_Tor_finite}
    Let $f \colon X \to S$ be a morphism of derived schemes over $k$.
    We let
\[
\DCoh(X / S) \subset \QCoh(X)
\]
be the full subcategory of almost perfect complexes with bounded Tor-dimension relative to $S$.
\end{definition}

\begin{notation}
    In the setting of \cref{notation:Coh_Tor_finite}, if $S= \Spec(A)$ is a derived affine scheme, we simply denote $\DCoh(X / \Spec(A))$ by $\DCoh(X / A)$.
\end{notation}

\begin{remark}\label{rem:relative_Coh_over_base/id}
    In the setting of \cref{notation:Coh_Tor_finite}, for $f= p \colon X \to \Spec(k)$ the structure morphism, we have
\[
\DCoh(X /  k) \simeq \DCoh(X) \ .
\]
For $f= id_X \colon X \to X$ we have
\[
\DCoh(X /  X) \simeq \Perf(X)
\]
by the characterization of perfect complexes as those almost perfect complexes with finite Tor-amplitude.
	In general, the $\infty$-category $\DCoh(X /  S)$ is an enhancement of Lieblich's relative perfect complexes \cite{Lie06}.
\end{remark}
 
\begin{definition}\label{def:proper_support_QCoh}
    Let $f \colon X \to S$ be a morphism of derived schemes locally almost of finite presentation over $k$, and let $F \in \QCoh(X)$.
    We say that $F$ has \emph{proper support relative to $S$} if there exists a closed subscheme $ Z \hookrightarrow X$ such that $Z \to S$ is proper and $F |_{X \setminus Z} \simeq 0$.
    We denote by $\QCoh_{prop/S}(X) \subset \QCoh(X)$ the stable full subcategory spanned by quasi-coherent complexes with proper support relative to $S$.
\end{definition}

    In the next definition, our main case of interest will be in the setting where $X$ is a scheme, $\Spec(A) \in \dAff_k$ is a derived affine scheme, and  $f \colon X_A \to X$, $g \colon X_A \to \Spec(A)$ are the projection morphisms.
\begin{definition}\label{def:proper_support_DCoh_rel}
    Let $f \colon X \to S$ and $g \colon X \to T$ be morphisms of derived schemes, locally almost of finite presentation over $k$.
    We denote by $\DCoh_{prop/T}(X/S)$ the pullback 
\[
    \DCoh_{prop/T}(X/S) \coloneqq \QCoh_{prop/T}(X) \times_{\QCoh(X)} \DCoh(X/S) \ .
\]
\end{definition}

\begin{notation}
    In the setting of \cref{def:proper_support_DCoh_rel}, let $f = id_X \colon X \to X$.
    We have $\DCoh(X/X) = \Perf(X)$ by \cref{rem:relative_Coh_over_base/id}, and in this case we simply write
\[
\Perf_{prop / T}(X) \coloneqq\DCoh_{prop/ T}(X/X) \ .
\]
\end{notation}

\begin{theorem}[{\cite[Theorem 5.0.5]{BZNP17}}]\label{thm:transforms_Coh}
	Let $X$ be a derived algebraic space, locally of finite presentation over $k$, and let $\Spec(A) \in \dAff_k$ locally of finite presentation over $k$.
	Then we have an equivalence
\begin{align*}
\DCoh_{prop/A}(X_A / & X) \stackrel{\sim}{\to} \Funst_k(\DCoh(X), \DCoh(\Spec(A))) \\
& F \mapsto p_{A, \ast} (F \otimes_{\OO_{X_A}} p_X^\ast (-)) \ ,
\end{align*}
where $p_X \colon X_A \to X$ and $p_A \colon X_A \to \Spec(A)$ are the projection morphisms.
\end{theorem}

\begin{corollary}\label{coro:points_IndCoh}
	Let $X$ be a derived algebraic space locally of finite presentation over $k$, and let $\Spec(A) \in \dAff_k$.
	If $A$ is a classical regular ring of finite presentation over $k$, then
\[
\IndCoh(X)_A^{pp} \simeq \DCoh_{prop / A}(X_A /  X) \ .
\]
\end{corollary}

\begin{proof}
	Since $\IndCoh(X)$ is auto-dual (see \cref{recollection:auto_duality_IndCoh}), we have an equivalence
\[
\IndCoh(X)_A \simeq \IndCoh(X)^\vee \otimes_k \Mod_A \simeq \Funst(\DCoh(X), \Mod(A)) \ ,
\]
inducing an equivalence $\IndCoh(X)_A^{pp} \simeq \Funst_k(\DCoh(X), \Perf(A))$.
	Since $A$ is regular noetherian, we have $\Perf(A) \simeq \DCoh(\Spec(A))$.
	Hence, the claim follows by \cref{thm:transforms_Coh}.
\end{proof}

\begin{lemma}\label{lemma:perfect_implies_realtively_perfect}
	Let $f \colon X \to Y$ be a morphism of derived schemes over $k$.
	Assume that $X$ is noetherian and eventually coconnective, and $f$ has finite Tor-amplitude.
    Then, we have an inclusion
\[
\Perf(X) \subset \DCoh(X / Y) \  .
\]
    Moreover, if $g \colon X \to T$ is a morphism of derived schemes locally almost of finite presentation over $k$, the above restricts to  an inclusion $\Perf_{prop/T}(X) \subset \DCoh_{prop/T}(X / Y).$
\end{lemma}

\begin{proof}
	Since $X$ is noetherian and eventually coconnective, we have $\Perf(X) \subset \DCoh(X)$.
	Let $F \in \Perf(X)$.
	We need to prove that $F$ has finite relative Tor-amplitude relative to $Y$, that is, there exist integers $\alpha \leq \beta$ such that, for each $M \in \QCoh(Y)^\heartsuit$, we have $\pi_i(F \otimes_{\OO_X} f^\ast (M)) = 0$ for all $i \notin [\alpha, \beta]$.
	By assumption, there exist integers $a \leq b$ such that $F$ has absolute Tor-amplitude contained in $[a,b]$ and a positive integer $n$ such that $f$ has Tor-amplitude $\leq n$.
	In particular, for each $M \in \QCoh(Y)^\heartsuit$ there exists an integer $k \leq n$ such that $f^\ast (M) \in \QCoh(X)_{\leq k, \geq 0}$.
	We show by induction on $0 \leq k \leq n$ that $\pi_i(F \otimes_{\OO_X} f^\ast (M)) = 0$ for all $i \notin [a, b+k]$.
	This directly implies that $F$ has relative Tor-amplitude contained in $[a, b+n]$.\\ \indent
	If $k=0$, then $f^\ast (M) \in \QCoh(X)^\heartsuit$.
	Since $F$ has absolute Tor-amplitude $[a,b]$, we get that $\pi_i(F \otimes_{\OO_X} f^\ast (M)) = 0$ for all $i \notin [a, b]$.
	Assume now that $0 < k \leq n$ and that we have proven our claim for $k-1$.
	Consider the fiber/cofiber sequence
\begin{equation}\label{eq:truncation}
\tau_{\geq k} f^\ast (M) \to f^\ast (M) \to \tau_{< k} f^\ast (M) \ .
\end{equation}
	We have that $\tau_{\geq k} f^\ast (M) \in \QCoh(X)^\heartsuit[k]$ and $\tau_{\leq (k-1)} f^\ast (M) \in \QCoh(X)_{< k, \geq 0}$.
	By the case $k=0$, we have that $\pi_i(F \otimes_{\OO_X} \tau_{\geq k} f^\ast (M)) = 0$ for all $i \notin [a+k, b+k]$.
	Using the induction hypothesis, the claim follows by passing to the long exact sequence of homotopy groups in \eqref{eq:truncation}. The statement in the properly supported case is immediate.
\end{proof}

\begin{remark}\label{rem:Perf_to_Coh}
	Let $X$ be a derived scheme over $k$.
    If $X$ and $\Spec(A) \in \dAff_k$ are eventually coconnective and locally of finite presentation over $k$,
    then \cref{lemma:perfect_implies_realtively_perfect} implies that there is an inclusion 
\begin{align}
    \Perf(X_A) \subset \DCoh(X_A /  X) \ .
\end{align}
    Similarly, we have $\Perf_{prop/A}(X_A) \subset \DCoh_{prop/A}(X_A / X).$
\end{remark}

\begin{lemma}\label{lemma:pullback_relative_coherent}
	Consider a pullback square of derived schemes over $k$
\[
\begin{tikzcd}[sep=small]
	{Y} && {X} \\
	\\
	{T} && {S}
	\arrow["{g}", from=1-1, to=1-3]
	\arrow["{q}"', from=1-1, to=3-1]
	\arrow["p", from=1-3, to=3-3]
	\arrow["f"', from=3-1, to=3-3]
\end{tikzcd}
\]
	with $f$ (and $g$) affine morphisms.
	Then, the functor $g^\ast \colon \QCoh(X) \to \QCoh(Y)$ restricts to
\[
g^\ast \colon \DCoh(X/S) \to \DCoh(Y / T) \  .
\]
\end{lemma}

\begin{proof}
	The proof of \cite[Lemma B.4]{DPS+26} applies verbatim to our setting.
\end{proof}

\begin{lemma}\label{lemma:push_relative_perf}
	Let $f \colon X \to S$ be a separated morphism of derived schemes.
	If $f$ is locally almost of finite presentation, then the functor $f_\ast \colon \QCoh(X) \to \QCoh(S)$ restricts to a functor
\[
f_\ast \colon \DCoh_{prop/S}(X/S) \to \Perf(S) \ .
\]
\end{lemma}

\begin{proof}
    Assume first that $f$ is proper.
	Let $F \in \DCoh(X/S)$ be of relative Tor-amplitude contained in $[a,b]$, with $a \leq b$ integers, and $M \in \QCoh(S)^\heartsuit$.
	By \cite[Theorem 5.6.0.2]{Lur18}, it is enough to prove that $f_\ast(F)$ has finite Tor-amplitude.
	By projection formula \cite[Corollary 6.3.4.3]{Lur18}, we have an equivalence
\[
f_\ast(F) \otimes_{\OO_S} M \simeq f_\ast(F \otimes_{\OO_X} f^\ast M) \ .
\]
	Since $F \otimes_{\OO_X} f^\ast M \in \QCoh(X)_{\geq a, \leq b}$ and $f$ is cohomologically bounded by \cite[Theorem 3.4.2.1 \& Proposition 2.5.4.4]{Lur18}, the result follows.
	More precisely, $f_\ast(F)$ is perfect of Tor-amplitude contained in $[a-n, b]$ where $n \in \N$ is such that $f_\ast(\QCoh(X)_{\geq 0}) \subset \QCoh(S)_{\geq -n}$.\\ \indent
    We now show the general claim by reducing to the case of a proper morphism.
    By \cite[Chapter 5, Proposition 2.1.6]{GR17}, the morphism $X \to S$ factors as $X \xrightarrow{j} \overline{X} \to S$, with $j$ an open immersion and $\overline{X} \to S$ proper.
    We will show that $j_\ast \colon \QCoh(X) \to \QCoh(\overline{X})$ restricts to $j_\ast \colon \DCoh_{prop/S}(X/S) \to \DCoh(\overline{X}/S)$, so that we are reduced to the case of a proper morphism as desired.
    Consider $E \in \DCoh_{prop/S}(X/S)$. That $j_\ast(E)$ is in $\DCoh(\overline{X}/S)$ can be checked Zariski locally on $\overline{X}$.
    Let $Z \subset X$ be a closed subscheme such that $E |_{X \setminus Z} \simeq 0$.
    We choose the open cover $\{j(X), \widebar{X} \setminus j(Z)\}$ of $\overline{X}$.
    The restriction of $j_\ast E$ to both opens is relative coherent: for $j(X)$ it holds by hypothesis, and for $\widebar{X}\setminus j(Z)$ we have the equality $j_\ast(E) |_{\widebar{X}\setminus j(Z)} \simeq 0$ by flat base change.
\end{proof}

\begin{proposition}\label{prop:FM_DCoh}
	Consider the pullback square of derived schemes over $k$
\[
\begin{tikzcd}[sep=small]
	{Y} && {X} \\
	\\
	{T} && {S}
	\arrow["{g}", from=1-1, to=1-3]
	\arrow["{q}"', from=1-1, to=3-1]
	\arrow["p", from=1-3, to=3-3]
	\arrow["f"', from=3-1, to=3-3]
\end{tikzcd}
\]
	Assume that $p$ (hence, $q$) is separated and locally almost of finite presentation and $f$ (hence, $g$) is affine.
	If $F \in \Perf_{prop / T}(Y)$, then the Fourier--Mukai transform
\[
\Phi_F \colon \DCoh(X/S) \to \QCoh(T) \  .
\]
takes values in $\Perf(T)$.
\end{proposition}

\begin{proof}
	By \cite[Theorem 5.6.0.2]{Lur18}, we have that $\Phi_F$ takes values in $\APerf(T)$.
	Hence, we need to show that for fixed $G \in \DCoh(X/S)$, we have that $\Phi_F(G) \in \QCoh(T)$ has finite Tor-amplitude over $T$.
	By \cref{lemma:pullback_relative_coherent}, we have $g^\ast G \in \DCoh(Y/T)$.
	Since $F$ is perfect over $Y$, we have that $F \otimes_{\OO_Y} g^\ast G \in \DCoh(Y/T)$.
    Moreover, we have that $F \otimes_{\OO_Y} g^\ast G$ has proper support over $T$ since $F$ does.
    Indeed, if $U \subset Y$ is an open subscheme, complement of a closed subscheme of $Y$ proper over $T$ and such that $F|_U \simeq 0$, we have 
\[
(M \otimes F )|_U \simeq M|_U \otimes F|_U \simeq 0 \ ,
\]
	for each $M \in \DCoh(Y / T)$.
	The result then follows by \cref{lemma:push_relative_perf}.
\end{proof}

\begin{corollary}\label{coro:Perf_is_pseduo_perf}
	Let $X$ be a separated derived scheme locally almost of finite presentation over $k$ and $\Spec(A) \in \dAff_k$.
	If $F \in \Perf_{prop/A}(X_A)$, then the Fourier--Mukai transform
\[
\Phi_F \colon \DCoh(X) \to \Mod_A \  .
 \]
takes values in $\Perf(A)$.
\end{corollary}

\begin{proof}
	The morphism $f \colon \Spec(A) \to \Spec(k)$ is affine, $X$ is locally almost of finite presentation over $k$ by assumption, and $\DCoh(X / \Spec(k)) = \DCoh(X)$ by \cref{rem:relative_Coh_over_base/id}.
	Hence, the claim follows directly from \cref{prop:FM_DCoh}.
\end{proof}

\subsection{\texorpdfstring{Representability of $\mathbf{Perf}_{prop}(X)$}{Representability of Perf-prop(X)}}

    Let $X$ be a derived scheme.
    The goal of this section is to prove that the stack $\mathbf{Perf}_{prop}(X)$, see \cref{def:moduli_Perf_prop}, is locally geometric and compute its cotangent complex.
    The strategy is as follows: we choose a compactification $\overline{X}$ of $X$, and show that the natural morphism $\mathbf{Perf}_{prop}(X) \to \mathbf{Perf}_{prop}(\overline{X})$ is an open immersion.
    As the latter is locally geometric, this is enough to conclude.
    To start, we need to show that the construction of $\mathbf{Perf}_{prop}(-)$ is functorial with respect to open immersions.
\begin{lemma}\label{lem:pushforward_perfprop}
    Let $j \colon X \hookrightarrow Y$ be an open immersion of derived schemes over a scheme $S$.
    Assume that $Y$ is separated over $S$.
    The functor $j_\ast \colon \QCoh(X) \to \QCoh(Y)$ restricts to a functor $j_\ast \colon \Perf_{prop/S}(X) \to \Perf_{prop/S}(Y)$.
\end{lemma}

\begin{proof}
    Let $F \in \Perf_{prop/S}(X)$ be a perfect complex on $X$ topologically supported on a closed subscheme $Z \xhookrightarrow{i} X$, proper over $S$.
    Denote by $U$ the complement of $Z$ in $X$, and by $V$ the complement of $j(Z)$ in $Y$.
    Notice that the composition $j \circ i \colon Z \to Y$ is a monomorphism since both $i$ and $j$ are.
    Moreover, since properness depends only on the underlying truncation, the morphism $j \circ i$ is proper by \cite[\href{https://stacks.math.columbia.edu/tag/01W6}{Lemma 01W6}]{stacks-project}.
    Hence, by \cite[\href{https://stacks.math.columbia.edu/tag/04XV}{Lemma 04XV}]{stacks-project} we have that $Z$ is a closed subscheme of $Y$, proper over $S$.
    Thus, it is enough to show that $j_\ast F$ is perfect and $j_\ast(F)|_V \simeq 0$.
    That $j_\ast(F)|_V \simeq 0$ follows from flat base change applied to the diagram 
    \begin{align}
        \begin{tikzcd}[ampersand replacement = \&]
            U \ar{r} \ar[hook, d] \& V \ar[hook, d] \\
            X \ar[r, hook, "j"] \& Y \ .
        \end{tikzcd}
    \end{align}
    Hence, we have $j_\ast(F) \in \QCoh_{prop/S}(Y)$.
    To verify that $j_\ast(F)$ is perfect, it is sufficient to check this on a Zariski cover of $Y$. As a cover of $Y$, we take $\{j(X), V\}$: the restriction of $j_\ast(F)$ to $V$ is, by the above, equivalent to the zero sheaf, hence perfect.
    The restriction $j^\ast j_\ast (F)$ is also perfect: indeed, the counit morphism $j^\ast j_\ast(F) \to F$ is an equivalence as the pushforward functor $j_\ast$ is fully faithful.
    In particular, $j^\ast j_\ast(F)$ is perfect. 
\end{proof}

\begin{corollary}\label{coro:open_imm_Perf_prop}
    Let $j \colon X \to \overline{X}$ be an open immersion of derived schemes over $k$.
    Assume that $\overline{X}$ is separated.
    Then, the pushforward $j_\ast$ induces a morphism of derived stacks 
    \begin{align*}
        j_\ast \colon \mathbf{Perf}_{prop}(X) \to \mathbf{Perf}_{prop}(\overline{X}) \ .
    \end{align*}
\end{corollary}

\begin{proof}
    After noticing that $\overline{X}_A = \overline{X} \times_k \Spec(A)$ is separated, this follows directly by applying \cref{lem:pushforward_perfprop} to the open immersion $X_A \hookrightarrow \overline{X}_A$ of schemes over $\Spec(A)$, $\Spec(A) \in \dAff_k$.
\end{proof}

\begin{remark}
    Let $j \colon X \to Y$ be as in the statement of \cref{lem:pushforward_perfprop}. The proof of \cref{lem:pushforward_perfprop} shows that the Tor-amplitude of $j_\ast(F)$ as a sheaf on $Y$ is the same as the Tor-amplitude of $F$ as a sheaf on $X$. In particular, given two integer numbers $a \leq b$, if we denote by $\mathbf{Perf}_{prop}^{[a,b]}(X)$ the derived stack parametrizing perfect complexes with Tor-amplitude in $[a,b]$, then $j_\ast$ restricts to a morphism
    \[j_\ast \colon \mathbf{Perf}_{prop}^{[a,b]}(X) \to \mathbf{Perf}_{prop}^{[a,b]}(Y) \ .\]
\end{remark}

\begin{lemma}\label{lemma:representability_Perf_porper_case}
    Let $X$ be a proper derived scheme, locally almost of finite presentation over $k$.
    Then $\mathbf{Perf}(X)$ is locally geometric and locally of finite presentation over $k$.
    Moreover, for any morphism $x_F \colon \Spec(A) \to \mathbf{Perf}(X)$, corresponding to a perfect complex $F \in \Perf(X_A)$, we have an equivalence
\[
\mathbb{T}_{\mathbf{Perf}(X), F} \simeq \Hom_{QCoh(X_A)}(F,F)[1] \ .
\]
\end{lemma}

\begin{proof}
    The first part is \cite[Theorem 5.1.1]{HalpernLeistnerPreygel2023}. 
    Thanks to \cite[Lemma 2.3 \& Proposition 2.7]{Porta25_HKR}, one can compute explicitly the cotangent complex of $\mathbf{Perf}(X)$ at the point $x_F$: unraveling the definitions in \emph{loc. cit.} and dualizing, one gets the result. 
\end{proof}

\begin{lemma}\label{lemma:push_repr_open_immersion}
    Let $j \colon X \to \overline{X}$ be an open immersion of separated derived schemes locally almost of finite presentation over $k$.
    Then, the morphism $j_\ast \colon \mathbf{Perf}_{prop}(X) \to \mathbf{Perf}_{prop}(\overline{X})$
of \cref{coro:open_imm_Perf_prop} is representable by an open immersion.
\end{lemma}

\begin{proof}
    Let $x_F \colon \Spec(A) \to \mathbf{Perf}_{prop}(\overline{X})$ be a morphism in $\dSt_k$ corresponding to a perfect complex with proper support $F \in \Perf_{prop / A}(\overline{X}_A)$.
    Let $Z \subset \overline{X}_A$ be a proper derived scheme over $\Spec(A)$ such that $F_{\big|V} \simeq 0$, where $V \coloneqq \overline{X}_A \setminus Z$.
    Up to shrinking $Z$, we may assume that $V$ is the biggest open in $\overline{X}_A$ where $F$ vanishes.
    Indeed, since $F$ is perfect, we can consider $W$ the maximal open subscheme of $\overline{X}_A$ where $F$ vanishes.
    Setting $Z' \coloneqq (\overline{X}_A \setminus W)_{\text{red}}$, 
    we have that $Z'$ is a closed subscheme of $Z$.
    In particular, we have that $Z'$ is proper over $\Spec(A)$.
    The topological subspace $\big| Z \big| \cap (\big|\overline{X}_A \setminus X_A\big|)$ is closed in $\big|Z\big|$.
    As the projection morphism $p \colon \overline{X}_A \to \Spec(A)$ becomes proper when restricted to $Z$, the image $\big|p\big| \big(\big| Z \big| \cap \big|\overline{X}_A \setminus X_A\big|\big)$ is a closed subset of $\Spec(A)$.\\ \indent
    Set $U \coloneqq \Spec(A) \setminus \big|p\big| \big(\big| Z \big| \cap \big|\overline{X}_A \setminus X_A\big|\big)$.
    We claim that there is an equivalence
\[
U \simeq \mathbf{Perf}_{prop}(X) \times_{\mathbf{Perf}_{prop}(\overline{X})} \Spec(A) 
\]
    compatible with the open immersion $U \hookrightarrow \Spec(A)$ and such that the projection $U \to \mathbf{Perf}_{prop}(X)$ corresponds to $F_{\big| p^{-1}(U)}$, which is perfect and with proper support over $U$ as $F$ is perfect and with proper support over $A$.
    This is enough to conclude the proof.
    To do so, let $T \in \dAff_k$ be an affine derived scheme, and let $f \colon T \to \Spec(A)$ be a morphism. The perfect complex $F_T \coloneqq (\text{id}_X \times f)^\ast F \in \mathbf{Perf}_{prop}(\overline{X})(T)$ lies in $\mathbf{Perf}_{prop}(X)$ if and only if the support of $F_T$ lies in $X_A \times_{\Spec(A)} T$.
    We notice that the support is compatible with pullbacks, i.e., we have 
\[
F_T|_{V \times_{\Spec(A)}T} \simeq 0 
\]
    and $V \times_{\Spec(A)}T$ is the biggest open where $F_T$ vanishes.
    In particular, we see that the support of $F_T$ lies in $X \times T$ if and only if 
    \[
    Z \times_{\Spec(A)} T \subset \overline{X}_A \times_{\Spec(A)} T
    \]
    factors through the open $X_A \times_{\Spec(A)} T \subset \overline{X}_A \times_{\Spec(A)} T$, which is equivalent to the condition that $f$ factors through $U$.
\end{proof}

\begin{proposition}\label{repr_Perf_prop}
    Let $X$ be a derived scheme, separated and locally almost of finite presentation over $k$.
    Then $\mathbf{Perf}_{prop}(X)$ is locally geometric and locally of finite presentation over $k$.
    Moreover, for any morphism $x_F \colon \Spec(A) \to \mathbf{Perf}_{prop}(X)$, corresponding to a perfect complex with proper support $F \in \Perf_{prop /A}(X_A)$, we have an equivalence
\[
\mathbb{T}_{\mathbf{Perf}_{prop}(X), F} \simeq \Hom_{QCoh(X_A)}(F,F)[1] \ .
\]
\end{proposition}

\begin{proof}
    By \cite[Chapter 5, Proposition 2.1.6]{GR17}, we can factor the structure morphism $p \colon X \to \Spec(k)$ as
\[
X \xrightarrow{j} \overline{X} \xrightarrow{\overline{p}} \Spec(k)
\]
where $j$ is an open immersion and $\overline{p}$ is proper and locally almost of finite presentation.
   Since $
j_\ast \colon \mathbf{Perf}_{prop}(X) \to \mathbf{Perf}_{prop}(\overline{X}) \simeq \mathbf{Perf}(\overline{X})$ is representable by a Zariski open immersion by \cref{lemma:push_repr_open_immersion}, the first claim follows by \cref{lemma:representability_Perf_porper_case}.

    Let now $x_F \colon \Spec(A) \to \mathbf{Perf}(X)$, corresponding to a perfect complex with proper support $F \in \Perf_{prop /A}(X_A)$.
    Since open immersions are étale and since $j_{A, \ast} \colon \QCoh(X_A) \to \QCoh(\overline{X_A})$ is fully faithful, we have an equivalence
\[
\mathbb{T}_{\mathbf{Perf}_{prop}(X), F} \simeq \mathbb{T}_{\mathbf{Perf}(\overline{X}), j_{A,\ast} F} \simeq \Hom_{\QCoh(\overline{X}_A)}(j_{A,\ast} F, j_{A,\ast} F)[1] \simeq \Hom_{\QCoh(X_A)}(F, F)[1] \ ,
\]
where the second equivalence follows by \cref{lemma:representability_Perf_porper_case}.
\end{proof}

\subsection{\texorpdfstring{The equivalence $\mathbf{Perf}_{prop}(X) \simeq \M_{\IndCoh(X)}$}{The equivalence Perf-prop(X) = M-IndCoh(X)}}\label{Appendix:equivalence}

    In this paragraph we prove that there is an equivalence in $\dSt_k$
\[
i \colon \mathbf{Perf}_{prop}(X) \to \M_{\IndCoh(X)} \,
\]
for $X$ a classical separated scheme locally of finite presentation over $k$.
    More precisely, we construct an étale morphism $i$ as above, which is an open immersion and is surjective on $k$-points.

\begin{construction}\label{construction:Perf_to_Coh}
    Let $X$ be an eventually coconnective and separated derived scheme, locally almost of finite presentation over $k$ and $\Spec(A) \in \dAff_k$.
    If $\Spec(A)$ is eventually coconnective and locally of finite presentation over $k$, then $\Perf(A) \subset \DCoh(\Spec(A))$ and thus the $\infty$-groupoid $\M_{\IndCoh(X)}(\Spec(A))$ is identified with the full $\infty$-subgroupoid 
\[
\IndCoh(X)_A^{pp} \subset \Funst_k(\DCoh(X), DCoh(\Spec(A))) \simeq \DCoh_{prop/A}(X_A / X)
\]
    spanned by those $F \in \DCoh_{prop/A}(X_A / X)$ such that the associated Fourier--Mukai transform $\Phi_F \colon \DCoh(X) \to \DCoh(A)$ takes values in $\Perf(A)$, where the second equivalence is \cref{thm:transforms_Coh}.\\ \indent  
    By \cref{rem:Perf_to_Coh} we have an inclusion of full subcategories of $\QCoh(X_A)$
\[
i_A \colon \Perf_{prop/A}(X_A) \to \DCoh_{prop/A}(X_A / X) \ .
\]
    Moreover, by \cref{coro:Perf_is_pseduo_perf} the morphism $i_A$ factors through $\IndCoh(X)_A^{pp}$.
    Hence, we get a fully faithful functor
    \begin{align}\label{eq:morph_Perf_to_Coh}
        i_A \colon \Perf_{prop/A}(X_A)^\simeq \to \M_{\IndCoh(X)}(\Spec(A)), \quad F \mapsto \phi_F \ .
    \end{align}
    By \cref{repr_Perf_prop}, the stack $\mathbf{Perf}_{prop}(X)$ is nilcomplete (\cite[Definition 17.3.2.1]{Lur18}) and locally of finite presentation over $k$, i.e., determined by its values on eventually coconnective simplicial rings, locally of finite presentation over $k$.
    If the same holds for $\M_{\IndCoh(X)}$ (e.g., if  $\IndCoh(X)$ is of finite type), then \eqref{eq:morph_Perf_to_Coh} defines a morphism of derived stacks over $k$
\[
 i \colon \mathbf{Perf}_{prop}(X) \to \M_{\IndCoh(X)} \ .
\]
\end{construction}

\begin{proposition}\label{proposition:open_immersion}
    Let $X$ be an eventually coconnective and separated derived scheme, locally almost of finite presentation over $k$.
    Assume that $\IndCoh(X)$ is of finite type (e.g., if $X$ is classical by \cref{coro:moduli_IndCoh}).
    Then, the morphism $i \colon \mathbf{Perf}_{prop}(X) \to \M_{\IndCoh(X)}$ of \cref{construction:Perf_to_Coh} is an open immersion.
\end{proposition}

\begin{proof}
    Notice that $\mathbf{Perf}_{prop}(X)$ and $\M_{\IndCoh(X)}$ are locally of finite presentation over $k$ by \cref{repr_Perf_prop} and \cref{repr_moduli_of_objects}.
    Moreover, they are nilcomplete by \cite[Proposition 2.7.3.2-(c) \& Proposition 19.2.1.5]{Lur18} and the fact that having proper support is a topological property.
    Thus, the morphism of \cref{construction:Perf_to_Coh} is well defined and is determined by its value on eventually coconnective simplicial rings, locally of finite presentation over $k$.
    We will show that $i$ is étale and a monomorphism.
    Since
\[
i_A \colon \Perf_{prop/A}(X_A) \to \IndCoh(X)_A^{pp}
\]
    realizes $\Perf_{prop/A}(X_A)$ fully faithfully as a subcategory of $\IndCoh(X)_A^{pp}$ for $A$ eventually coconnective and locally of finite presentation over $k$ (\cref{construction:Perf_to_Coh}), it follows that the diagonal $\Delta_i$ of $i$ is an equivalence.
    Thus $i$ is a monomorphism.\\ \indent
    We now prove that $i$ is étale.
    Since both $\mathbf{Perf}_{prop}(X)$ and $\M_{\IndCoh(X)}$ are locally of finite presentation over $k$, it is enough to prove that $i$ is formally étale.
    Since $\mathbf{Perf}_{prop}(X)$ is locally of finite presentation, it is enough to show that the relative cotangent complex is trivial after pulling back to any derived affine locally of finite presentation.
    Since the cotangent complexes of $\mathbf{Perf}_{prop}(X)$ and $\M_{\IndCoh(X)}$ are perfect, we can equivalently show that the relative tangent complex is trivial.
    Moreover, by \cite[Corollary 2.7.4.4]{Lur18}, it is enough to prove that the relative tangent complex is zero after pulling back to any classical (hence almost coconnective) affine scheme, locally of finite presentation over $k$.
    Let then $F \in \Perf_{prop/A}(X_A)$ be an $A$-point of $\mathbf{Perf}_{prop}(X)$, with $A$ almost coconnective and locally of finite presentation over $k$.
    We have identifications
\[
\mathbb{T}_{\mathbf{Perf}_{prop}(X), F} \simeq \Hom_{\Perf(X_A)}(F,F)[1], \qquad \mathbb{T}_{\M_{\IndCoh(X)}, i(F)} \simeq \Hom_{\DCoh(X_A/X)}(F,F)[1] \ ,
\]
and both $\Hom$'s can be computed in $\QCoh(X_A)$.
    Indeed, we have inclusions of full subcategories $\Perf(X_A) \subset \QCoh(X_A)$, and $\IndCoh(X)_A^{pp} \subset \DCoh(X_A/X) \subset \QCoh(X_A)$ by \cref{rem:Perf_to_Coh}.
    Since the morphism $i(A)$ is fully faithful, it follows that the tangent morphism
\[\mathbb{T}(i) \colon \mathbb{T}_{\mathbf{Perf}_{prop}(X), F} \to \mathbb{T}_{\M_{\IndCoh(X)}, i(F)}
\]
is an equivalence.
    Hence, the morphism $i$ is formally étale.
\end{proof}

    In order to prove that $i$ is an equivalence, we will need the following result.
    It is well-known, but we could not locate a reference and we include a proof here for completeness.
    Recall that, throughout all the paper, the field $k$ is assumed to be algebraically closed.

\begin{lemma}\label{lemma:open_and_surj_is_iso}
    Let $i \colon \mathcal{U} \to \X$ be an open immersion of locally geometric derived Artin stacks that are locally almost of finite presentation over $k$.
    If $i$ is surjective on $k$-points, then $i$ is an equivalence.
\end{lemma}

\begin{proof}
    Since $i$ is an open immersion, the derived structures of $\mathcal{U}$ and $\X$ agree.
    Hence, it is enough to prove that an open immersion of classical locally geometric higher stacks locally of finite presentation over $k$, which is surjective on $k$-points, is an equivalence.
    For this, since $\mathcal{U}$ and $\X$ are assumed to be locally of finite presentation over $k$, it is enough to prove that for each $\Spec(A) \to \X$, with $A$ classical and locally of finite presentation over $k$, the pullback $i_A \colon \mathcal{U} \times_\X \Spec(A) \to \Spec(A)$ is an equivalence.
    Notice that $i_A$ is an open immersion of classical schemes, which is surjective on $k$-points, since open immersions of stacks are schematic and the properties "open immersion" and "surjective on $k$-points" are stable under pullbacks.
    To prove that $i_A$ is an equivalence, we show that the closed subscheme $Z \coloneqq \Spec(A) \setminus (\mathcal{U} \times_\X \Spec(A)) \to \Spec(A)$, say with its reduced structure, is empty.
    Indeed, if $Z$ is non-empty, then it must admit a $k$-point since $k$ is algebraically closed and $Z$ is locally of finite presentation over $k$.
    Our assumption that $i$ is surjective on $k$-points implies that $Z(k)$ is empty, hence that $\mathcal{U} \simeq \X$.
\end{proof}

\begin{theorem}\label{thm:comparison_moduli_Perf_DCoh}
    Let $X$ be a separated, eventually coconnective derived scheme, locally of finite presentation over $k$.
    Assume that $\IndCoh(X)$ is of finite type (e.g., if $X$ is classical by \cref{thm:IndCoh_finite}).
    Then, the morphism $i \colon \mathbf{Perf}_{prop}(X) \to \M_{\IndCoh(X)}$ is an equivalence.
\end{theorem}

\begin{proof}
    By \cref{proposition:open_immersion}, we have that $i$ is an open immersion.
    By \cref{rem:relative_Coh_over_base/id} and \cref{coro:points_IndCoh}, we also have that $i$ induces an equivalence on $k$-points.
    The statement then follows by \cref{lemma:open_and_surj_is_iso}.
\end{proof}

\section{\texorpdfstring{A fully commutative proof of \cref{coro:Lagr_strong}}{A fully commutative proof of coro:Lagr-strong}}\label{appendix:the_case_of_Perf}

	In this paragraph, we give a proof of \cref{coro:Lagr_strong}, from which \cref{thm:main_DCoh} follows, but directly from the point of view of shifted Lagrangian morphisms of derived stacks.
    That is, we make no reference to relative Calabi--Yau structures.
    The proof is very similar, and in most places we only need to give references to \cite{Cal15} instead of \cite{BD19, BD21}.
	We will show the analogue of \cref{coro:Lagr_strong}, and then proceed to apply \cref{thm:composition_lagr_morph} as in the proof of \cref{thm:main_DCoh}.
	
\begin{definition}[{\cite[Definition 2.1]{PTVV13}}]\label{def:O_compact}
	Let $\X \in \dSt_k$.
	We say that $\X$ is $\OO$-compact if, for any $\Spec(A) \in \dAff_k$, we have:
\begin{enumerate}\itemsep=0.2cm
	\item $\OO_{\X_A}$ is a compact object in $\QCoh(\X_A)$;
	\item For any $F \in \Perf(X_A)$, we have $\Gamma(\X_A, F) \in \Perf(A)$.
\end{enumerate}
\end{definition}

\begin{recollection}
	Let $\X \in \dSt_k$ be $\OO$-compact and $[\X] \colon \Gamma(\X, \OO_\X) \to k[-d]$ be a morphism in $\Mod_k$, $d \in \ZZ$.
	Then for any $\Spec(A) \in \dAff_k$ and $F \in \Perf(\X_A)$, we can consider the composite map
\[
\Gamma(\X_A, F) \otimes_A \Gamma(\X_A, F^\vee) \to \Gamma(\X_A, \OO_{\X_A}) \xrightarrow{[\X] \otimes_k id_A} A[-d] \ .
\]
	By $\OO$-compactness of $\X$, we have $\Gamma(\X_A, F^\vee) \in \Perf(A)$, so that the above composite gives rise to a map
\begin{equation}\label{eq:cup_orientation}
- \cap [\X] \colon \Gamma(\X_A, F) \to \Gamma(\X_A, F^\vee)^\vee[-d] \ .
\end{equation}
\end{recollection}

	In the setting of derived algebraic geometry, the following notion of boundary condition is the replacement of relative left Calabi--Yau structures that we will need:
	
\begin{definition}[{\cite[Definition 2.4]{PTVV13}}]\label{def:boundary_condition}
	Let $\X \in \dSt_k$ be an $\OO$-compact derived stack over $k$ and let $[\X] \colon \Gamma(\X, \OO_\X) \to k[-d]$ be a morphism in $\Mod_k$, $d \in \ZZ$.
	We say that $[\X]$ is an $\OO$-orientation of degree $d$ if the map $- \cap [\X] \colon \Gamma(\X_A, F) \to \Gamma(\X_A, F^\vee)^\vee [-d]$ in \eqref{eq:cup_orientation} is an equivalence for any $\Spec(A) \in \dAff_k$ and $F \in \Perf(\X_A)$.
\end{definition}

	Orientations give rise to shifted symplectic structure on mapping stacks via AKSZ construction: 

\begin{theorem}[{\cite[Theorem 2.5]{PTVV13}}]\label{thm:AKSZ}
	Let $\Y \in \dSt_k$ be a geometric stack equipped with an $n$-shifted symplectic form $\omega$.
	Let $\X \in \dSt_k$ be an $\OO$-compact derived stack equipped with an $\OO$-orientation $[X] : \Gamma(\X, \OO_\X) \to k[-d]$ of degree $d$.
	Assume that the derived mapping stack $\mathbf{Map}(\X, \Y) \in \dSt_k$ is itself geometric and locally of finite presentation over $k$. 
	Then, $\mathbf{Map}(\X, \Y)$ carries a canonical $(n-d)$-shifted symplectic structure.
\end{theorem}

\begin{remark}
	In the setting of \cref{thm:AKSZ}, let $\Spec(A) \in \dAff_k$. Consider a point $x \colon \Spec(A) \to \mathbf{Map}(\X, \Y)$, corresponding to a morphism $f \colon \X_A \to \Y$.
	Then, the underlying non-degenerate pairing of \cref{thm:AKSZ} is described in the proof of \cite[Theorem 2.5]{PTVV13}, and has the following description on the tangent space of $\mathbf{Map}(\X, \Y)$ at the point $x$:
\[
\Gamma(\X_A, f^\ast \mathbb{T}_\Y) \wedge \Gamma(\X_A, f^\ast \mathbb{T}_\Y) \xrightarrow{f^\ast(\omega)} \Gamma(\X_A, \OO_{X_A}) [n] \simeq \Gamma(\X, \OO_{X}) \otimes_k A [n] \xrightarrow{[\X] \otimes_k id_A} A[n-d] \ .
\]
\end{remark}

\begin{remark}\label{rem:AKSZ_locally_geometric}
	As noticed in the introduction of \cite[Section 2.3]{PTVV13}, \cref{thm:AKSZ} remains valid if we replace geometric stacks with locally geometric stacks.
	Hence, we might take $\Y = \mathbf{Perf}$ in the statement of \cref{thm:AKSZ}.
	See \cite[Theorem 2.12 \& Corollary 3.13]{PTVV13}.
\end{remark}

\begin{lemma}\label{lem:O-compact}
	Let $X$ be a proper Gorenstein scheme over $k$.
	Then $X$ is $\OO$-compact.
\end{lemma}

\begin{proof}
	We need to verify items $(1)$ and $(2)$ of \cref{def:O_compact}.
	Item $(1)$ follows for example from \cite[Proposition 9.6.1.1]{Lur18}.
	For item $(2)$, we need to show that the pushforward along the projection $p_A \colon X_A \to \Spec(A)$ preserves perfect complexes.
	By \cite[Proposition 9.6.1.1]{Lur18}, this is equivalent to showing that $p_{A, \ast} \colon \QCoh(X_A) \to \Mod_A$ preserves compact objects.
	Since $p \colon X \to \Spec(k)$ is proper and (flat, hence) of finite Tor-amplitude, the same holds for its basechange $p_A$.
	Thus, the result follows by combining \cite[Theorem 5.6.0.2, Proposition 6.1.3.1]{Lur18} with \cite[Proposition 7.2.4.23-(4)]{Lur17}.
\end{proof}

	Absolute Calabi--Yau structures give rise to symplectic forms (\cref{thm:CY_lagr}), and to construct Lagrangian structures one needs a relative notion.
	Similarly, orientations allow us to construct shifted symplectic forms (\cref{thm:AKSZ}), and to construct Lagrangian structures one needs the following relative notion introduced by Calaque \cite{Cal15}.
	
\begin{definition}\label{def:boundary_structure}
	Let $\X, \Z \in \dSt_k$ be $\OO$-compact and $f \colon \X \to \Z$.
	Let $[\X] \colon \Gamma(\X, \OO_\X) \to k[-d]$ be an orientation.
	A boundary structure on $f$ with respect to $[\X]$ is a homotopy $h$ to $0$ of the composite map in $\Mod_k$
\[
f_\ast [\X] \colon \Gamma(\Z, \OO_\Z) \to \Gamma(\X, \OO_\X) \xrightarrow{[\X]} k[-d] \ .
\]
\end{definition}

\begin{recollection}
	In the setting of \cref{def:boundary_structure}, let $F \in \Perf(\Z)$.
	Consider the pairing in $\Mod_k$
\begin{equation}\label{eq:pairing}
\Gamma(\Z, F) \otimes \Gamma(\X, f^\ast F^\vee) \to \Gamma(\X, f^\ast F) \otimes \Gamma(\X, f^\ast F^\vee) \to \Gamma(\X, \OO_\X) \xrightarrow{[\X]} k[-d] \ .
\end{equation}

	Since $f^\ast F^\vee \in \Perf(\X)$, then \eqref{eq:pairing} gives rise to a map
\[
g \colon \Gamma(\X, f^\ast F^\vee) \to \Gamma(\Z, F)^\vee[-d] \ ,
\]
and we denote by $\Gamma(\Z, F)^\bot$ the fiber of $g$.
	Consider then the commutative diagram in $\Mod_k$
\[
\begin{tikzcd}[sep=small]
	{\Gamma(\Z, F)^\bot} &&&& \\
	\\
	{\Gamma(\Z, F^\vee)} && {\Gamma(\X, f^\ast F^\vee)} \\
	\\
	&& {\Gamma(\X, f^\ast F)^\vee[-d]} && {\Gamma(\Z, F)^\vee[-d]}
	\arrow[from=1-1, to=3-3]
	\arrow["{\delta_h}", dashed, from=3-1, to=1-1]
	\arrow["can", from=3-1, to=3-3]
	\arrow[from=3-3, to=5-3]
	\arrow["g", from=3-3, to=5-5]
	\arrow[from=5-3, to=5-5]
\end{tikzcd}
\]
	By definition the homotopy $h$ gives rise to a homotopy to $0$ of the composite $g \circ can$, hence defining the dotted morphism
\begin{equation}\label{eq:boundary_morphism}
\delta_h \colon \Gamma(\Z, F^\vee) \to \Gamma(\Z, F)^\bot \ .
\end{equation}
\end{recollection}

\begin{definition}\label{def:non_deg_boundary}
	Let $\X, \Z \in \dSt_k$ be $\OO$-compact and $f \colon \X \to \Z$.
	Let $[\X] \colon \Gamma(\X, \OO_\X) \to k[-d]$ be an orientation and $h$ be a boundary structure for $f$ with respect to $[\X]$.
	We say that $h$ is non-degenerate if, for any $F \in \Perf(\Z)$, the morphism $\delta_h \colon \Gamma(\Z, F^\vee) \to \Gamma(\Z, F)^\bot$ in \eqref{eq:boundary_morphism} is an equivalence.
\end{definition}

\begin{proposition}\label{prop:Lagr_strong_Perf}
	Let $X$ be a proper Gorenstein scheme over $k$ of dimension $d > 0$.
	Then, trivializations of $\omega_X^\bullet[-d]$ are in bijection with orientations of $X$
\[
[X] \colon \Gamma(X, \OO_X) \to k[-d] \ .
\]
    In this case, $p \colon X \to \Spec(k)$ carries a boundary structure $h$.
	Moreover, if $X$ is \textit{strong} Calabi--Yau, such an $h$ is non-degenerate.
\end{proposition}

\begin{proof}
	Since $\Gamma(X, \OO_X) \in \Perf(k)$ by \cref{lem:O-compact}, we have that giving a morphism $\Gamma(X, \OO_X) \to k[-d]$ is equivalent to giving a morphism
\[
k \to \Gamma(X, \OO_X)^\vee[-d] \simeq \Gamma(X, \omega^\bullet_X[-d]) \ .
\]
	Hence, morphisms $f \colon \Gamma(X, \OO_X) \to k[-d]$ are in bijection with $s_f \in \mathrm{H}^0\Gamma(X, \omega^\bullet_X[-d])$, which is easily seen to be an orientation if and only if $s_f$ is everywhere non-zero, that is, if $s_f$ induces an equivalence $\OO_X \simeq \omega_X^\bullet[-d]$.
	The proof of the Lagrangianity is carried out as in the proof of the relative left Calabi--Yau part in \cref{prop:CY_strong}.
\end{proof}

	In order to prove the analogue of \cref{thm:main_DCoh} for $\mathbf{Perf}(X) \coloneqq \mathbf{Map}(X, \mathbf{Perf})$, we need a variant of \cref{prop:Lagrangian_ev_fib}.

\begin{proposition}\label{prop:Lagrangian_ev_fib_Perf}
	Let $X$ be a Calabi--Yau scheme of dimension $d>0$.
	Then, the map
\[
(\fib, \ev_0, \ev_1) \colon \mathbf{Perf}_{prop}(X)^{\Delta^1} \to \mathbf{Perf}_{prop}(X)^{\times 3}
\]
is $(1-d)$-shifted Lagrangian.	
\end{proposition}

\begin{proof}
	This follows by \cite[Theorem 3.1.6]{Bu25} applied to $\mathfrak{X} = \mathbf{Perf}_{prop}(X)$, as 
\[
\mathbf{Perf}_{prop}(X)^{\Delta^1} \subset \mathrm{Filt}(\mathbf{Perf}_{prop}(X)) \coloneqq \mathbf{Map}([A^1 / \G_m], \mathbf{Perf}_{prop}(X))
\]
is an open substack.
\end{proof}

	We are now able to deduce the main result of this paragraph.

\begin{theorem}\label{thm:main_Perf}
	Let $X$ be a strong Calabi--Yau variety of dimension $d > 0$ and $i \colon \mathbf{\mathcal{L}} \to \mathbf{Perf}(X)$ be a $(1-d)$-shifted Lagrangian morphism.
	If $X$ is strong proper Calabi--Yau, then the morphism of derived stacks
\[
\fib \circ q \colon \widetilde{\mathbf{Syz}}^{\mathbf{\mathcal{L}}}(X) \to \mathbf{Perf}(X)
\]
is $(1-d)$-shifted Lagrangian.
\end{theorem}

\begin{proof}
	It follows directly by \cref{prop:Lagr_strong_Perf}, \cref{prop:Lagrangian_ev_fib_Perf} and \cref{thm:composition_lagr_morph}.
\end{proof}

\bibliographystyle{alpha}
\bibliography{bibliography}

\end{document}